\documentclass[reqno,a4paper]{amsart}
\usepackage{amssymb,stmaryrd}
\usepackage{amsfonts}
\usepackage{amstext}
\usepackage{algorithmic}
\usepackage{algorithm}
\usepackage{color}
\usepackage{graphicx}
\usepackage{epstopdf}
\usepackage[all]{xy}
\usepackage{MnSymbol}
\numberwithin{equation}{section}
\newtheorem{theorem}{Theorem}[section]
\newtheorem{lemma}[theorem]{Lemma}
\newtheorem{proposition}[theorem]{Proposition}
\newtheorem{corollary}[theorem]{Corollary}

\theoremstyle{definition}

\newtheorem{example}[theorem]{Example}

\theoremstyle{remark}
\newtheorem{remark}[theorem]{\bf{Remark}}

\newcommand{\Hom}{{\rm{Hom}}}

\newcommand{\R}{{\mathbb{R}}}

\newcommand{\C}{{\mathbb{C}}}
\newcommand{\Z}{{\mathbb{Z}}}

\newcommand{\<}{{\langle}}

\renewcommand{\>}{{\rangle}}

\newcommand{\CC}{{\mathcal{C}}}

\newcommand{\CR}{{\mathcal{R}}}

\newcommand{\CM}{{\mathcal{M}}}
\newcommand{\wedgeq}{{\wedge\kern-5pt\cdot}}

\newcommand{\cg}{{\mathfrak{g}}}

\renewcommand{\ker}{{\rm{ker}}}

\newcommand{\tens}{\otimes}

\newcommand{\id}{{\rm id}}

\newcommand{\extd}{{\rm d}}
\newcommand{\del}{{\partial}}
\newcommand{\eps}{\epsilon}
\newcommand{\ev}{{\rm ev}}

\newcommand{\und}{\underline}

\renewcommand{\imath}{\mathrm{i}}

\newcommand{\blu}[1]{{\color{blue}{#1}}}

\allowdisplaybreaks
\begin{document}

\title{Complex structure on quantum-braided planes}
\keywords{noncommutative geometry, quantum Riemannian geometry}

\subjclass[2020]{Primary 58B32; 81R50, 83C65, 46L87}

\author{Edwin Beggs}\author{Shahn Majid}\thanks{{\it Authors to whom correspondence should be addressed:} e.j.beggs@swansea.ac.uk  and s.majid@qmul.ac.uk}
\address{Department of Mathematics, Bay Campus, Swansea University, SA1 8EN, UK}
\address{Queen Mary University of London, 
School of Mathematical Sciences, Mile End Rd, London E1 4NS, UK}

\thanks{Ver 1}

\begin{abstract} We construct a  quantum Dolbeault double complex $\oplus_{p,q}\Omega^{p,q}$ on the quantum plane $\C_q^2$. This solves the long-standing problem that the standard differential calculus on the quantum plane is not a $*$-calculus, by embedding it as the holomorphic part of a $*$-calculus. We show in general that any Nichols-Woronowicz  algebra or braided plane $B_+(V)$, where $V$ is an object in an abelian $\C$-linear braided bar category of real type is a quantum complex space in this sense with factorisable Dolbeault double complex. We combine the Chern construction on $\Omega^{1,0}$ in such a Dolbeault complex for an algebra $A$ with its conjugate to construct a canonical metric compatible connection on $\Omega^1$ associated to a class of quantum metrics, and apply this to the quantum plane.  We also apply this to finite groups $G$ with Cayley graph generators split into two halves related by inversion, constructing such a  Dolbeault  complex $\Omega(G)$ in this case, recovering the quantum Levi-Civita connection for any edge-symmetric metric on the integer lattice with $\Omega(\Z)$ now viewed as a quantum complex structure. We also show how to build natural quantum metrics on $\Omega^{1,0}$ and $\Omega^{0,1}$ separately where the inner product in the case of the quantum plane, in order to descend to $\tens_A$,  is taken with values in an $A$-bimodule.  \end{abstract}
\maketitle 

\section{Introduction}

A long-standing puzzle which we address in this paper is what should be the natural $*$-differential calculus $(\Omega,\extd)$ on the quantum plane $\C_q^2$ with generators $x,y$ and relations $yx=qxy$ with $q$ real, and similarly for higher dimensional quantum planes. In noncommutative geometry, both from operator algebras\cite{Con} and from a more $*$-algebraic  `quantum Riemannian geometry'\cite{BegMa} approach, one effectively complexifies the analogue of smooth functions on a manifold and remembers the `real form' by the $*$-involution. If the algebra was commutative then the self-adjoint elements would recover the analogue of the real-valued functions. This is no problem for quantum planes and there are some natural such `real forms', but extending this to differential forms as a $*$-exterior algebra (so that $*$ commutes with the exterior derivative) has not been possible to achieve in a convincing manner. Particularly, the standard $\C_q[SL_2]$ or $\C_q[GL_2]$-covariant  2-dimensional calculus on $\C_q^2$ does not appear to be a $*$-calculus for an appropriate choice of $*$-structures on the different objects. 

Our solution to the problem is motivated by a similar issue for the lattice line $\Z$ where, as for any discrete space, differential 1-forms in noncommutative geometry are spanned as a vector space by the arrows of a bidirected graph, with $*$ naturally reversing arrows. This means that the smallest translation-invariant $*$-calculus $(\Omega(\Z),\extd)$ has an $\Omega^1$ that is 2 and not 1-dimensional over the algebra, which in turn allows for curvature and quantum gravity on these graphs. Moreover, $(\Omega(\Z),\extd)$ factorises into two sub-exterior algebras connected by $*$ forming a  factorisable Dolbeault double complex. The latter is a natural formulation in noncommutative geometry of a  `quantum complex manifold' in an algebraic sense\cite{BegSmi,BegMa}. We find exactly analogous results for the quantum plane $\C_q^2$, constructing a natural 4-dimensional $*$-calculus on it where the exterior algebra is a complex structure in this sense. The standard 2D calculus now appears as $\Omega^{\bullet,0}$ and, like  for the lattice line, this holomorphic (similarly the antiholomorphic) part is connected (i.e. the kernel of $\del$  is the constant functions). This is  different from classical complex geometry on the complex plane, but in other ways the structure is analogous. 

In fact, the quantum plane  has the additional structure of a braided plane or Hopf algebra in a braided category with primitive generators\cite{Ma},  leading to braided-derivatives and providing the key to our construction. Specifically, it can be cast in the form  $B_+(V)$, the braided-symmetric algebra or braided-plane (a.k.a. Nichols-Woronowicz algebra)  associated to  any rigid object $V$ in an Abelian braided category\cite{Ma:dbos}.  This has a standard differential calculus  constructed in  \cite{Ma:hod}\cite[Chap. 2]{BegMa}, which appears now as $\Omega^{0,\bullet}=B_+(V)\underline\tens \Lambda(V)$, where we use a braided-tensor product with the relevant exterior algebra  $\Lambda(V)$ on $V$. Our new result (which also works over any field if we are not interested in $*$) is  a doubled up version   
\[ \Omega=\Lambda(V)\underline\tens B_+(V)\underline\tens \Lambda(V)\]
 as a braided super tensor product. Working over $\C$, we require that $V$ is a unitary $*$-object in a braided bar category, using concepts from \cite{BegMa} and our earlier works, which then leads to a factorisable Dolbeault complex $\Omega(B_+(V))$ for all such braided planes, see Theorem~\ref{thmBV}. 

We then turn to an application of noncommutative  complex structures (in this Dolbeault complex sense) on an algebra $A$ that are factorisable at least to degree 2, namely to the construction of metric compatible connections. It is already known in this case\cite[Chap.7]{BegMa} that a hermitian metric on $\Omega^{1,0}$ leads to a  holomorphic `Chern connection' 
and we show in Section~\ref{secdoublech} how to combine this with its conjugate to obtain a hermitian metric compatible connection on $\Omega^1=\Omega^{1,0}\oplus\Omega^{0,1}$ for a class of generalised hermitian metrics on this. If this connection is $*$-preserving and torsion-free then we obtain a quantum Levi-Civita connection (QLC) on a corresponding ordinary (non-hermitian) quantum metric $(\ ,\ ): \Omega^1\tens\Omega^1\to A$ as needed for quantum Riemannian geometry. This is a powerful general result which does not suppose that $(\ ,\ )$ descends to $\tens_A$ as normally required, i.e. applies more generally.  Example~\ref{exZ} shows how this works for $\Omega(\Z)$ and we indeed recover  exactly the previously know unique QLC for the integer lattice in \cite{Ma:haw}. We then apply the same method in Section~\ref{secexplane} to the  quantum-braided plane $\C_q^2$ to obtain a natural hermitian metric compatible connection on $\Omega^1$ given by $\nabla=0$ on the image of $\del,\bar\del$, which can then be viewed as a QLC for a corresponding generalised quantum metric $(\ ,\ )$. This metric has a cross-diagonal form that is nonzero only on $\Omega^{1,0}\tens\Omega^{0,1}, \Omega^{0,1}\tens\Omega^{1,0}$, similarly to the $\Omega(\Z)$ case, but unlike that case does not descend to $\tens_A$. 

We are also able to adapt some of the braided category/diagrammatic methods used in the construction of hermitian metrics on $B_+(V)$ to now construct an ordinary (not hermitian) metric  $(\ , \ )$ on $\Omega^{1,0}$ alone. This works over any field. We can do the same on $\Omega^{0,1}$ and also double up to a metric on $\Omega^1$ that is nonzero only on $\Omega^{1,0}\tens\Omega^{1,0}, \Omega^{0,1}\tens\Omega^{0,1}$ and which, for the quantum plane, again does not descend to $\tens_A$. This time, however, we show for the quantum plane and similar examples how the quantum metric on $\Omega^{1,0}$ can be adapted so that it {\em does} descend to 
\[ (\ ,\ ):\Omega^{1,0}\tens_A\Omega^{1,0}\to L_+\]
 but now valued in a certain bimodule $L_+$ over $A$ rather than in $A$ itself. 
Similarly for $\Omega^{0,1}$ with $(\ ,\ )$ valued in a bimodule $L_-$. The two halves then combine to $(\ ,\ )$ on $\Omega^1$ to be valued in a bimodule $L=L_+\oplus L_-$ and now descends to $\tens_A$. The same $\nabla$ as before provides a metric compatible connection for this in a suitable sense adapted to the bimodule values. 

We note that complex structures in the Dolbeault double complex sense have been used in several places, such as \cite{Ma:spi},\cite{Hec},\cite{Bua}, while quantum Riemannian geometries in our sense with a quantum metric and bimodule connections have been used to construct models of quantum gravity, such as \cite{Ma:squ}. Bimodule connections themselves appeared in \cite{DV1,DV2,Mou}, while interesting works by others using bimodule connections include \cite{Car} and \cite{Sit} and, in slightly different formalisms \cite{Lan,AscWeb}.  Also, this style of noncommutative geometry is somewhat different from the operator algebras and spectral triples approach of Connes\cite{Con} but can be connected by a `geometric realisation' formalism \cite{BegMa:spe,Ma:can} whereby Dirac operators are built from spinor bundles as $A$-bimodules equipped with a bimodule connection, Clifford and $*$ structures. Our work brings the quantum plane and other braided planes into this active paradigm. It is also the case that the strong invertibility of a quantum metric in \cite{BegMa} imply that the metric element is central, and there is some interest in generalising to avoiding this.  We have demonstrated the point of view in \cite{BegMa} that a natural way to weaken the notion of quantum metric is to work with the inner product  form $(\ ,\ )$ and not require that it descends to $\tens_A$. We have then provided an novel alternative which is for the map to descend but to be valued in a bimodule or `line bundle'. which is also of interest classically\cite{BicKnu}. Another issue is that one has to solve for the QLC and this may not exist or be unique, leading to various proposals for variant frameworks, such as  \cite{Lan}. In our case, the Chern construction coming from a factorisable complex structure can indeed provide a natural metric-compatible connection, but not necessarily torsion-free.  It should also be mentioned that we have not imposed some form of symmetry on $(\ ,\ )$, and indeed our examples on the quantum plane q-deform something antisymmetric rather than symmetric.

An outline of the paper is as follows. In Section~\ref{secpre}, we recall some basics from quantum groups and differential structures with $*$. In particular, Lemma~\ref{lemqplane*} explains how to view the quantum plane as generated by $z,z^*$ with an action of a central extension of $\C_q[SU_{1,1}]$. We also recall the known construction of the standard 2D calculus as generalised to any $B_+(V)$ in \cite{Ma:hod} and denoted now as $(\Omega^{0,\bullet},\bar\del)$.  
Section~\ref{secdolb} then proceeds to the full double complex  $\Omega$ for such  braided planes obtained from unitary $*$-objects in a braided bar category.  Sections~\ref{secinner} and~\ref{secch} recap inner products and hermitian inner products and the Chern construction \cite[Chap.~7]{BegMa} on a factorisable double complex. A new result here is Proposition~\ref{propGdolb} which constructs a factorisable complex structure $\Omega(G)$ on a finite group where the Cayley graph suitably splits into disjoint arrows labeled by $\CC^{1,0},\CC^{0,1}\subset G$ connected by group inversion. 
We provide the order 12 group $A_4$ as a nontrivial example. Section~\ref{secdoublech} then provides a `doubled up' Chern construction for bimodule connections $\nabla$ on $\Omega^1=\Omega^{1,0}\oplus\Omega^{0,1}$ that are  compatible with a class of quantum metrics. Section~\ref{secbundle} covers the case where inner products are valued in a bimodule $L$. This section has some remarks for the hermitian setting, but the main focus is on constructing examples in the non-hermitian $(\ ,\ )$ case.

\section{Preliminaries}\label{secpre}

We will work in a Abelian braided tensor category but the examples of interest will be the $k$-linear case of comodules over a quasitriangular Hopf algbra $H$. Here $k$ is the ground field and we will be interested in the case $k=\C$. Without going into formal details, there is a tensor product of objects $\tens$ and for any two objects $V,W$ a braiding natural isomorphism $\Psi_{V,W}:V\tens W\to W\tens V$ which we will denote diagrammatically by $\Psi=\includegraphics{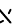}$ reading down the page. We also denote $\Psi^{-1}=\includegraphics{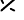}$ and by the coherence theorem we can work as if the category is strictly associative (inserting brackets as needed). This underlies the diagrammatic notation\cite{Ma:alg} where $\tens$ is denoted by juxtaposition, an algebra $B$ in the braided category has product denoted $\includegraphics{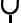}$ and a bialgebra has this as well as a coproduct denoted $\includegraphics{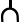}$. There is an associated unit $\eta:\underline{1}\to B$ and a counit $\eps:B\to \underline{1}$ where the unit object $\underline{1}$ (in our case the ground field) is denoted by omission. Given $B,C$ two algebras in the category there is a braided tensor product $B\underline{\tens}C$ as in \cite{Ma:alg} which is best described diagrammatically but in the concrete case where objects are sets,  appears as $(b\tens c).(b'\tens c')= b\Psi_{C,B}(c\tens b')c'$ for all $b,b'\in B$ and $c,c'\in C$. For a bialgebra we need $\eps,\Delta$ to be algebra homs where $\Delta: B\to B\underline{\tens}B$. There is also a notion of a braided-antipode, see \cite{Ma:alg}. 

If $V$ is an object in a braided category, we define braided-factorials $[n,\Psi]!: V^{\tens n}\to V^{\tens n}$ as introduced in \cite{Ma} and 
\[ B=B_+(V)={TV\over \oplus_n \ker[n,\Psi]!}\]
where $TV$ is the tensor algebra on $V$ in the category as the associated `braided symmetric algebra'. There are different routes to this, but the one in \cite{Ma:dbos} was to start with $TV$ as a braided Hopf algebra with $V$ primitive (so $\Delta v=v\tens 1+1\tens v$ for $v\in V$ in the concrete case).  Then if $V$ is rigid then $TV$ and $TV^*$ are dually paired as braided-Hopf algebras via $[n,\Psi]!$, so $B_+(V)$ is the cokernel on one side, with a similar (dual) braided-Hopf algebra on the other side. The evaluation between $TV$ and $TV^*$ is given by
\[
\ev_V{}^{(n)}\big( [n,\Psi]! \tens\id^n\big) : V^{\tens n} \tens V^*{}^{\tens n}\to \underline{1},
\]
where $\ev_V:V\tens V^*\to  \underline{1}$ and $\ev_V{}^{(n)}$ denotes nesting $n$ of these with no crossings. The algebras $B_+(V)$ and a similar $B_-(V)$ (where we use $-\Psi$) are also called `Nichols-Woronowicz' algebras (in an equivalent form where the braided factorials are written as a sum over the symmetric group). However, the braided category setting and the braided-Hopf algebra structure was not present in any of their works, being introduced in \cite{Ma:dbos}. The construction also applies to construct $U_q(n_-)$ where it recovers a braided-Hopf algebra implicit (without a general theory of these) in the work of Lusztig. We consider $B=B_+(V)$ as a noncommutative `coordinate algebra' and define $\varpi: B\to V$ to be the projection to degree 1 and
\[ \bar\del=(\id\tens\varpi)\Delta: B\to B\tens V.\]
These are the braided differentials as in \cite{Ma} but in a `coaction' form. In \cite{Ma:hod}\cite[Prop.~2.78]{BegMa} it explained that this defines a differential calculus on $\Omega^{0,1}=B\tens V$ with differential $\bar\del$ except that we did not use an (anti)holomorphic notation as here. This extends to higher degree to give a DGA with degrees 
\[\Omega^{0,\bullet}= B\underline{\tens} \Lambda(V),\quad \Lambda(V)={TV\over \<{\rm Im}(\id+\Psi)\>}    \]
where $ \<{\rm Im}(\id+\Psi)\> $ denotes the two sided ideal generated by $ {\rm Im}(\id+\Psi) $ and 
where the latter is the braided exterior algebra that we need. We denote its product by $\wedge$. The exterior derivative here is defined similarly by
\[ \bar\del=(\id\tens\wedge)(\id\tens\varpi\tens\id)(\Delta\tens \id): B\underline{\tens} \Lambda^q\to B\underline{\tens}\Lambda^{q+1}.\]
which makes sense as elements of $\Omega^{0,1}$ of the form $1\tens V$ are exact since $\del v=1\tens v$ for $v\in V\subset B$ given the primitive form of the coproduct. Hence $\bar\del|_\Lambda=0$. 

We will also need the notion of a bar category\cite{BegMa:bar}\cite[Sec. 2.8]{BegMa}. This is a monoidal or in our case braided monoidal category with the further structure of a `bar functor' that expresses complex conjugation. This comes with a natural equivalence $\Upsilon: \bar{\ }\circ \tens\Rightarrow \tens^{op}\circ (\bar{\ }\times\bar{\ })$ obeying a hexagon identity, and some supplementary data. To keep things concrete, we focus on the relevant case of $\CM^H$ the comodules over a  $*$-Hopf algebra $H$. Recall that the latter is a Hopf algebra equipped with an antilinear antimultiplicatve $*$-involution (making it a $*$-algbra) which commutes with the coproduct $\Delta:H\to H\to H$ and $\eps:H\to \C$. In this case $(S\circ *)^2=\id$. Then $\CM^H$ this becomes a bar category as follows. Given $V\in \CM^H$, we first define $\overline{V}$ as the same set but a vector space with conjugate action of $\C$, 
\[ \overline{\lambda v}=\bar\lambda \, \overline{v}\]
for all $v\in V$ and $\lambda\in \C$ with $\bar\lambda$ its complex conjugate. Here $\overline{V}$ denotes $v\in V$ viewed in $\overline{V}$. As a comodule,
\begin{equation}\label{conjcoact} \Delta_{\overline{V}}\overline{V}= \overline{v_0}\tens v_1^*\end{equation}
and the natural equivalence is 
\[ \Upsilon_{V,W}(\overline{v\tens w})=\overline{w}\tens\overline{v}\]
for all $v\in V, w\in W$. If $\phi:V\to W$ is a comodule map then $\bar\phi:\overline{V}\to\overline{W}$ is $\bar \phi(\overline{v})=\overline{\phi(v)}$ for all $v\in V$.

Next, recall that a  coquasitriangular Hopf algebra means one equipped with a convolution-invertible linear map $\CR: H\tens H\to \C$ (in our case) subject to the axioms
\[ h_1 g_1 \CR(g_2, h_2)=\CR(g_1,h_1) g_2 h_2\]
\[ \CR(fg,h)=\CR(f,h_1)\CR(g,h_2),\quad \CR(f,gh)=\CR(f_1,h)\CR(f_2,g)\]
for all $f,g,h\in H$, where $\Delta h=h_1\tens h_2$ (sum of such terms understood). In the $*$-Hopf algebra case $\CR$ is of real type\cite{Ma} if
\[ \overline{\CR(g,h)}=\CR(h^*,g^*)\]
for all $h,g\in H$. The category of right comodules $\CM^H$ is braided via $\Psi_{V,W}(v\tens w)=w_0\tens v_0\CR(v_1\tens w_1)$, 
where $\Delta_V v=v_0\tens v_1$ (sum of terms understood) is the right coaction and in the real type case one has
\[ \Upsilon_{V,W}\circ \overline{\Psi_{V,W}}=\Psi_{\overline{V},\overline{W}}\circ\Upsilon_{W,V}.\]
Finally, a $*$-object in a bar category means an object $V$ equipped with a morphism $\star:V\to \overline{V}$ such that $\bar\star\circ\star={\rm bb}_V$ where ${\rm bb}_V:V\to \overline{\overline{V}}$ is a natural transformation that is part of the data of a bar category but which or $\CM^H$ is just the identity under the identification of the underlying sets. This comes down to $V$ a comodule equipped with an antilinear involution $*$ such that $(*\tens\id) \Delta_V(v)= \Delta_{\overline{V}}(v^*)$ which given (\ref{conjcoact}), just amounts to
\[ (*\tens *)\Delta_V v= \Delta_V v^*\]
for all $v\in V$ when the underlying sets are identified (so that we do not write the overlines explicitly here).  This is familiar in the case where $V$ is a $*$-algebra which is a comodule algebra, and says that $H$ coacts {\em unitarily}; the notion makes sense just for comodules equipped with $*$ and we say that $V$ is then a {\em unitary comodule}. If $V,W$ are both unitary then the above property of the braiding translates to
\begin{equation}\label{psidag} \Psi_{V,W}\circ\dagger=\dagger\circ\Psi_{W,V},\quad \dagger={\rm flip}(*\tens *)\end{equation}
for the relevant $*$ on $V,W$. The flip map here is the content of $\Upsilon$ after identifications. We have translated the notion of a $*$-object and the properties we need to an ordinary algebraic language in the case of $\CM^H$. Henceforth, we will use this ordinary language with the underlying bar category taking care of the categorical book-keeping behind the scenes. Using this, it is known\cite[Prop.~2.110]{BegMa} that $B_+(V)$ (when $B$ is unitary) is a *-braided-Hopf algebra. The axioms for this in our ordinary language are that $B$ has an involution whereby 
\[ *\circ \cdot=\cdot\circ \dagger\]
(i.e. antimultiplicative), which just amounts to $B_+(V)$ a $*$-algebra, and at the same time
\[ \dagger\circ\Delta=\Delta \circ *\]
(i.e. anticomultiplicative). There is also a condition on the braided-antipode which holds. The reason for this result is evident in degree 2 where the product of $B_+(V)$ quotients by the kernel of $\id-\Psi$. Hence, since $\dagger$ on $TV$ commutes with  this, it descends to $B_+(V)$ to define $*$ as expected in the product. The coproduct on degree 2 is likewise given by $\Delta(vw)=vw\tens 1+ 1\tens vw + (\id+\Psi)(v\tens w)$ and we again use (\ref{psidag}). For higher degrees we have to use properties of the braided-factorials. 

We next recall that for $q\ne 0$ real, $\C_q[SU_{1,1}]$ denotes $\C_q[SL_2]$ with its standard generators $a,b,c,d$ and standard coquasitriangular structure\cite{Ma} as a $*$-Hopf algebra of real type with 
\[  a^*=d, \quad d^*=a, \quad b^*=q^{-1} c,\quad c^*=q b.\]
We also let $\C_{q^\alpha}[\R_{>0}]=\C[t,t^{-1}]$ as usual with group-like generator $t$ and coquasitriangular structure $\CR(t^m,t^n)=q^{\alpha mn}$ but now with $t^*=t$ and $\alpha$ real. This is also a $*$-Hopf algebra of real type. 

\begin{lemma}\label{lemqplane*}  Let $H=\C_q[SU_{1,1}]\tens\C_{q^\alpha}[\R_{>0}]$ with $q>0$. Then $V=\C^2$ with basis $z,w$ is a unitary $H$-comodule via
\[ (z\ w)\mapsto (z\ w)\tens\begin{pmatrix}a & b\\ c & d\end{pmatrix}\tens t,\quad z^*=q^{{1\over 2}}w,\quad w^*=q^{-{1\over 2}}z\]
and $\C_q^2$ as generated by $z,z^*$ with relations $z^* z=q z z^*$ is a $*$-algebra in the bar category $\CM^H$. Moreover, for $\alpha=3/2$, this is a $*$-braided Hopf algebra in  $\CM^H$ as a braided bar-category of real type.\end{lemma}
\proof  For the first part
\[ (*\tens *)\Delta_V z=(*\tens *)((z\tens a+w\tens c)t)=q^{{1\over 2}}( w\tens d+z\tens b)t=q^{{1\over 2}}\Delta_V w=\Delta_V z^*,\]
\[ (*\tens *)\Delta_V w=(*\tens *)((w\tens d+z\tens b)t)=q^{-{1\over 2}}( z\tens a+w\tens c)t=q^{-{1\over 2}}\Delta_V z=\Delta_V w^*.\]
This automatically extends to a $*$-algebra with a unitary coaction of $H$ (i.e. the coaction is a $*$-algebra map).
Next, it is already known that with $\alpha=3/2$ the braiding $\Psi$ coming from the $R$ matrix of $\C_q[SL_2]$ with the `quantum group normalisation' times an extra factor $q^{3/2}$  from $\CR$ for $\C_q[\R_{>0}]$ becomes 
\[\Psi_{V,V}(z\tens z)= q^2 z\tens z,\quad \Psi_{V\tens V}(w\tens w)=q^2 w\tens w,\quad \Psi_{V\tens V}(z\tens w)=qw\tens z,\]
\[ \Psi_{V,V}(w\tens z)=qz\tens w + (q^2-1)w\tens z\]
which gives $B_+(V)=\C_q^2$ as a braided-Hopf algebra. From the theory above, since $V$ is unitary, this is a $*$-braided Hopf algebra also. \endproof

It is also known\cite{Ma:hod}\cite[Chap.~2]{BegMa} that the standard 2D-calculus on the quantum plane is recovered by the above construction, which we will note denote $\Omega^{0,\bullet}$, namely with basis $\bar\del z,\bar\del w$ and relations
\[ (\bar\del z) z=q^2 z\bar\del z,\quad (\bar\del w) w=q^2 w\bar\del w,\quad (\bar\del z) w=q w\bar\del z,\]
\[ (\bar\del w) z=q z\bar\del w+ (q^2-1)w\bar\del z,\quad (\bar\del z)^2=(\bar\del w)^2=0,\quad \bar\del w\wedge\bar\del z+q^{-1}\bar\del z\wedge\bar\del w=0.\]
This is not, however, a $*$-calculus unless $q=1$.

\section{Factorisable Dolbeault complex $\Omega^{p,q}$ on braided planes}\label{secdolb}

Continuing from the description of $\Omega^{0,1}$ in \cite{Ma:hod}\cite{BegMa}, we now make the right handed version using $\del$ and $\Omega^{1,0}=V\tens B$, where $B=B_+(V)$ as before. Similarly to the above, we have an exterior algebra 
\[ \Omega^{\bullet,0}=\Lambda(V)\underline{\tens} B, \quad \del=(-1)^p(\wedge\tens\id)(\id\tens\varpi\tens\id)(\id\tens\Delta): \Lambda^p\underline{\tens} B\to \Lambda^{p+1}\underline{\tens} B\]
where the $(-1)^p$ is needed to convert the right-handed construction to  left graded derivation. Apart from the sign, the proofs are those in \cite{Ma:hod}\cite{BegMa} with the diagrams reflected in a mirror and all braid-crossings then reversed. We now combine these constructions.

\begin{figure}
\[ \includegraphics[scale=0.75]{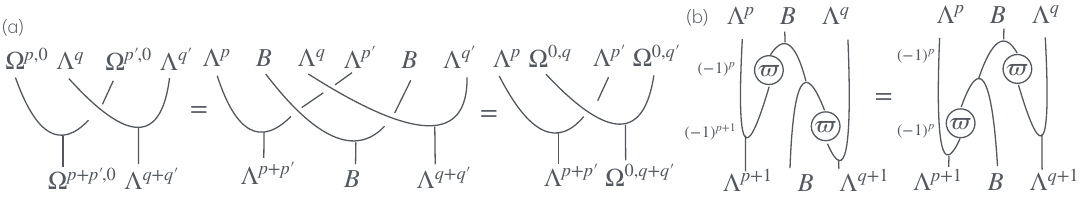}\]
\caption{Proof of  Proposition~\ref{propdoublecomp}. Here $B$ stands for $B_+(V)$. \label{fig1}}
\end{figure}

\begin{proposition}\label{propdoublecomp} Let $V$ be an object in an Abelian braided category and $B=B_+(V)$, $\Lambda=\Lambda(V)$, $(\Omega^{0,q},\bar\del)$ and $(\Omega^{p,0},\del)$ DGAs as above. Then 
\[ \Omega^{p,q}=\Lambda^p\underline{\tens} B\underline{\tens} \Lambda^q,\]
\[ \del:\Omega^{p,q}\to \Omega^{p+1,q},\quad \del=(-1)^p((\wedge\tens \id)(\id\tens\varpi\tens\id)(\id\tens \Delta))\tens\id,\]
 \[\bar\del: \Omega^{p,q}\to \Omega^{p,q+1},\quad \bar\del=(-1)^p\id\tens((\id\tens\wedge)(\id\tens\varpi\tens\id)(\Delta\tens\id))\]
 make $\Omega^{p,q}$ into a double complex.\end{proposition}
\proof The diagrams behind the construction are shown in Figure~\ref{fig1}. In (a) we use associativity of the braided tensor product to group the product of $\Omega^{p,q}=\Omega^{p,0}\und\tens\Lambda^q=\Lambda^p\und\tens\Omega^{0,q}$ with $\Omega^{p',q'}$, in two ways. From the first form, at left, it is clear that $\del$ extends as a derivation. We use the derivation property on $\Omega^{p+p',0}$ and for the 2nd term of the result we use that  $\del$ is a morphism to push it to the top. Similarly for $\bar\del$ a derivation. Part (b) then shows that $\del,\bar\del$ anti-commute, using coassociativity of the coproduct. The bimodule action of $b\in B$ is as $1\tens b\tens 1$ using the braided tensor product algebra structure. \endproof

\begin{corollary}\label{corfac}  In the context of Proposition~\ref{propdoublecomp}, $\wedge:\Omega^{p,0}\tens_B \Omega^{0,q}\to \Omega^{p,q}$ is a bimodule isomorphism by the canonical identification of $B\tens_B B=B$. Moreover, $\wedge:\Omega^{0,q}\tens_B\Omega^{p,0}\to \Omega^{p,q}$ is a also a bimodule isomorphism.
\end{corollary}
\proof  The first is automatic as $(\Lambda^p\tens B)\tens_B(B\tens\Lambda^q)=\Lambda^p\tens B\tens\Lambda^q$ and was already used in the proof of Proposition~\ref{propdoublecomp}.  We next consider
the bimodule isomorphism
\[ \Theta^{p,q}:=(\Psi_{B\tens\Lambda^p}\tens \Psi_{\Lambda^q\tens B})(\id\tens\Psi_{\Lambda^q,\Lambda^p}\tens\id):\Omega^{0,q}\tens\Omega^{p,0}\to   \Omega^{p,0}\tens \Omega^{0,q},\]
 before we take $\tens_B$, where
\begin{equation}\label{psirev} \Psi_{\Lambda^q,B}: \Omega^{q,0}\to \Omega^{0,q},\quad \Psi_{B,\Lambda^p}:\Omega^{0,p}\to \Omega^{p,0}\end{equation}
are respectively a right $B$-module isomorphism (as in the proof above) and left $B$-module isomorphism.  This is such that $\wedge\Theta^{p,q}$ is exactly the second wedge computed using the braided tensor product structure, which also descends to $\tens_B$. Next, $\Psi_{\Lambda^q,B}^{-1}$ is also a right module map and applied to the first factor gives an isomorphism 
\[ (B\tens\Lambda^q)\tens_B(\Lambda^p\tens B)\cong(\Lambda^q\tens B)\tens_B(\Lambda^p)\tens B=\Lambda^q\tens \Lambda^p\tens B\]
 by the standard identification of $B\tens_B$, so these vector spaces are isomorphic. Moreover, since $\Psi_{\Lambda^q,B}|_{\Lambda^q\tens 1}$ is trivial, we can identify $\Lambda^q\tens \Lambda^p\tens B$ with the subspace $1\tens\Lambda^q\tens \Lambda^p\tens B$ as representatives of $(B\tens\Lambda^q)\tens_B(\Lambda^p\tens B)$. It remains to note that
\[ \wedge\circ \Theta^{p,q}|_{1\tens \Lambda^q\tens \Lambda^p\tens B}=\Psi_{\Lambda^q,\Lambda^p\tens B},\]
which is invertible, so in this form  $\Psi_{\Lambda^q,\Lambda^p\tens B}^{-1}$ provides the inverse to the second wedge product. \endproof \endproof


This property of $\wedge$ being an isomorphism is called 
 {\em factorisability} in \cite[Chap.~7]{BegMa},  and in this case we denote inverse maps by $\Theta^{p00q}:\Omega^{p,q} \to \Omega^{p,0}\tens_B\Omega^{0,q}$ and 
$\Theta^{0qp0}:\Omega^{p,q} \to \Omega^{0,q}\tens_B\Omega^{p,0}$ as in \cite[Def.~7.21]{BegMa}. So far, our results work over any field. We now work over $\C$ and require $B=B_+(V)$ to be a $*$-algebra in a suitable sense. This needs the braided category to be a bar category as recalled in Section~\ref{secpre}. To keep things simple, we work as recalled there in the bar category of comodules over a coquasitriangular Hopf $*$-algebra $H$ of real type. We require $V$ to be a star-object in the bar category, which here just amounts to a unitary comodule. 

\begin{figure}
\[ \includegraphics[scale=0.9]{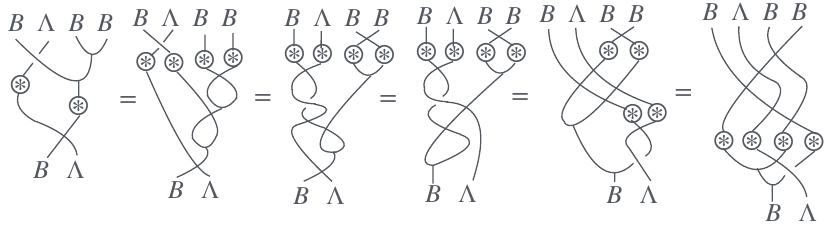}\]
\caption{Proof of  Proposition~\ref{thmBV}. Here $B$ stands for $B_+(V)$. \label{fig2}}
\end{figure}

\begin{theorem}\label{thmBV} If $(V,*)$ is a unitary right $H$-comodule then $B_+(V)$ is a $*$-algebra and $\Omega=\oplus_{p,q}\Omega^{p,q}$ and $\extd =\del+\bar\del$ is a $*$-differential calculus over $B_+(V)$. Moreover, this is the factorisable Dolbeault double complex  of a quantum complex structure in the sense of \cite[Chap. 7]{BegMa}.
\end{theorem}
\proof We already saw in Section~\ref{secpre} that $B=B_+(V)$ is a $*$-braided Hopf algebra. Both $V$ and $B$ are unitary and hence 
\[\dagger\circ\del =\dagger(\varpi\tens\id)\Delta=(\id\tens\varpi)\dagger\Delta=(\id\tens\varpi)\Delta\circ *=\bar\del\circ *\]
acting on $B$. This makes $\Omega^1=\Omega^{1,0}\oplus\Omega^{0,1}$ with $\extd=\del+\bar\del$ into a first order $*$-differential calculus provided the natural identification $(\Omega^{1,0})^*=\Omega^{0,1}$ is a skew-bimodule map, where we define $*$ on $\Lambda^1\underline\tens B$ by $\dagger$. This extends easily higher forms provided $(\Omega^{p,0})^*=\Omega^{0,p}$, where $\Lambda(V)$ also gets a well-defined graded-$*$ operation by $(v\wedge w)^*=-\wedge\dagger(v\tens w)$. The proof that we have a skew bimodule map is in Figure~\ref{fig2} where we show that
\[ b.(\omega\tens c).d=d^*.(\omega\tens c)^\dagger. b^*\]
for all $\omega\in \Lambda$, $b,c,d\in B$. We denote the flip map or vector space transposition by a solid crossing $\times$. The left hand side here is the first expression in the figure applied to $b\tens \omega\tens c\tens d$, since $b=(1\tens b)$ acts by the braided tensor product while $d$ acts by right multiplication. We use that $B$ is a $*$-algebra for the 2nd expression and (\ref{psidag}) for the 3rd. We the rearrange the operations that do not involve the braiding as these take place in the usual category of vector spaces for the 4th expression. For the 5th expression we likewise move the last two strands to before the braid crossing (regarding it as some linear map) and further rearrange the vector space operations to obtain the last expression, which corresponds to the above right hand side that we wanted to obtain. Thus, the two calculi are related by $*$. 

We similarly define $\Omega^m=\oplus_{p+q=m}\Omega^{p,q}$, $\extd=\del+\bar\del$ as part of the double complex and it is clear that this is then a $*$-DGA, where we define $*$ by iterated $\dagger$ so that
\[ (\omega\tens b\tens \eta)^*=(-1)^{|\omega||\eta|}\eta^*\tens b^*\tens \omega^*.\]
Moreover, by the remark after \cite[Prop.~7.5]{BegMa}, a complex structure is equivalent to having a double complex which is a $*$-calculus such that $\Omega^{p,q}{}^*=\Omega^{q,p}$ as a skew-bimodule map and $*\del=\bar\del *$. The proofs are similar to the degree 1 case covered.  Factorisability of the double complex is by Corollary~\ref{corfac}. \endproof

\begin{corollary}\label{corqplane}For $H=\C_q[SU_{1,1}]\tens \C_{q^{3\over 2}}[t,t^{-1}]$ as in Lemma~\ref{lemqplane*}, the quantum plane $\C_q^2$ generated by $z,z^*$ with $z^*z=qzz^*$ is a quantum complex manifold with Dolbeault double complex $\Omega^{\bullet,\bullet}$ with basis $\del z,\del z^*,\bar\del z,\bar\del z^*$ and relations in $\Omega^1$
\[ z\, \del z =q^2 (\del z)z,\quad  z^* \del z^*=q^2 (\del z^*)z^*,\quad  z\, \del z^* =q (\del z^*)z,\]
\[ z^*\del z=q (\del z)z^*+ (q^2-1)(\del z^*)z,\quad z\,\bar\del z^*=q^{-1} (\bar\del z^*)z+(q^{-2}-1)(\bar\del z)z^*,
 \]
 \[ z^* \bar\del z^* =q^{-2} (\bar\del z^*)z^*,\quad  z\, \bar\del z=q^{-2} (\bar\del z) z,\quad  z^* \bar\del z =q^{-1} (\bar\del z)z^*,\]
\[ (\bar\del z^*) z=q \, z\,\bar\del z^*+ (q^2-1)z^*\, \bar\del z,\quad (\del z)^2=(\del z^*)^2=0,\quad \del z^*\wedge\del z+q^{-1}\del z\wedge\del z^*=0\]
\[  \bar\del z \wedge \del z =-q^2 \del z\wedge \bar\del z,\quad  \bar\del z^* \wedge \del z^*=-q^2 \del z^*\wedge\bar\del z^*,\]
\[   \bar\del z\wedge \del z^* =-q \del z^*\wedge \bar\del z,\quad \bar\del z^*\wedge \del z=-q \del z\wedge\bar\del z^*- (q^2-1)\del z^*\wedge\bar\del z\]  
where $\bar\del z^*=(\del z)^*$ and the total derivative is $\extd=\del+\bar\del$. The factorisation structure is given by the braiding, i.e.
\[ \Theta^{1,1}(\bar\del z\tens\del z)=q^2 \del z\tens\bar\del z,\quad \Theta^{1,1}(\bar\del z^*\tens\del z^*)=q^2 \del z^*\tens\bar\del z^*,\quad \Theta^{1,1}(\bar\del z\tens\del z^*)= q\del z^*\tens\bar\del z\]
\[ \Theta^{1,1}(\bar\del z^*\tens\del z)=q\del z\tens\bar\del z^*+(q^2-1)\del z^*\tens\bar\del z.\]
\end{corollary}
\proof This now follows immediately by the general theory. Here $\Omega^{\bullet,0}$ is generated by   $\del z,\del w$ with relations  by similar calculations as for $\Omega^{0,\bullet}$ in \cite{Ma:hod},\cite[Chap.~2]{BegMa} and stated after Lemma~\ref{lemqplane*}. On the generators $1\tens z=\bar\del z, z\tens 1=\del z$ etc, the map $\Theta^{1,1}$ in Corollary~\ref{corfac} reduces to the braiding $\Psi$. We also changed to $z^*=q^{{1\over 2}}w$ (but this amounts to replacing $w$ by $z^*$ since all expressions are homogeneous in the degree of $w$), to obtain the relations stated. By the above proposition, this forms a $*$-calculus with cross relations as stated, again first computing with $w$ and then replacing by $z^*$. Moreover, by the proposition, we have a double complex with the required properties.  \endproof

\section{Chern connection on $\Omega^{1,0}$ and its doubling to $\Omega^1$}\label{secchern}

Here we look at some elements on quantum geometry made possible by the existence of a factorisable Dolbeault complex such as we constructed in  Corollary~\ref{corfac} and Theorem~\ref{thmBV}. For the rest of the paper, our focus is on noncommutative geometry and we will denote the coordinate algebra by $A$. The algebra $B=B_+(V)$  will be an example, but not the only one of interest. 

\subsection{Inner products on $\Omega^{1,0}$ and $\Omega^{0,1}$}\label{secinner}

We first consider the construction of hermitian inner products under the assumption that there is a $A$-bimodule map $\Theta^{1,1}$ as below (which for case of $B_+(V)$ we saw in the proof of factorisability
in Corollary~\ref{corfac}).  We then have the following commutative diagram for the solid arrows, in which all solid arrows are $A$-bimodule maps:
\begin{align} \label{eq47}
		\xymatrix{
\Omega^{1,0}\tens\overline{\Omega^{1,0} }\ \  \ar[r]^{\id\tens\medstar^{-1}}_{\cong}   &\ \  \Omega^{1,0}\tens\Omega^{0,1}    
 \ar[rd]^\wedge \ar@/^/[rrd]^{\phi_+}    &
			\\
	&		& \Omega^{1,1} \ar@{.>}[r]^{\phi\quad} 
	  & A      
	\\
	\Omega^{0,1}\tens\overline{\Omega^{0,1} }\ \  \ar[r]^{\id\tens\medstar^{-1}}_{\cong}
	  &\ \  \Omega^{0,1}\tens\Omega^{1,0}     \ar[ru]_\wedge    \ar[uu]^{\Theta^{1,1}}
 \ar@/_/[rru]_{\phi_-}  
}
\end{align}  
Here $\medstar^{-1}:\overline{\Omega^{1,0} } \to \Omega^{0,1} $ is given by $\medstar^{-1}(\overline{\xi})=\xi^*$, and we suppose that
$\Theta^{1,1}$ is invertible. This invertibility implies that $\phi_+$ is defined if and only $\phi_-$ is defined. The condition for e.g.\ $\phi_-$ to give a hermitian inner product $\<\xi,\overline{\eta}\>=\phi_-(\xi\tens\eta^*)$ is $\phi_-\circ\dagger=*\,\phi_-$, where $\dagger$ is the antilinear map from (\ref{psidag}). If we have $\dagger\circ \Theta^{1,1}=\Theta^{1,1}\circ\dagger$ then $\phi_+$ gives a hermitian inner product if and only $\phi_-$ does. One can also consider 
 the inner products $\phi_\pm$ to be valued in a $A$-bimodule $L$ with a star operation. 

 There is no reason for $\phi_+$ or $\phi_-$ to descend to an inner product over $\tens_A$, but if one of them does then so does the other, and we have the dotted line $\phi$ giving a commuting diagram. To see this, suppose that $\phi_-$ factors through 
$\Omega^{0,1}\tens_A\Omega^{1,0}$. By factorisability, this is the same as saying that it factors through $\Omega^{1,1}$, so we have the dotted map $\phi$ on the diagram, and then $\phi_+$ also factors through $\Omega^{1,1}$, and thus $\Omega^{1,0}\tens_A\Omega^{0,1}$. Moreover, if $\phi_\pm$ descend to $\tens_A$ then they give hermitian inner products if and only if $\phi$ is antihermitian. To see this, for
 $\xi,\eta$ in $\Omega^{0,1}$ or $\Omega^{1,0}$,
\[
\<\xi,\overline{\eta}\>^*=\phi(\xi\wedge\eta^*)^*=-\phi\big((\xi\wedge\eta^*)^*\big) = \phi(\eta\wedge\xi^*)=\<\eta,\overline{\xi}\>^*.
\]
Note that we do not require positivity. In fact, there is no obvious link between the positivity of the two inner products above, as one requires that $\phi(\xi\wedge\xi^*)\ge 0$ for all $\xi\in \Omega^{1,0}$ and the other that $\phi(\eta\wedge\eta^*)\ge 0$ for all $\eta\in \Omega^{0,1}$.

We can form a hermitian inner product on $\Omega^1$ by adding the inner products due to $\phi_+$ and $\phi_-$ and supposing that $\Omega^{1,0}$ and 
$\Omega^{0,1}$ are perpendicular. To be specific,  we use the direct sum $\Omega^1=\Omega^{1,0} \oplus\Omega^{0,1}$ to give projections $\pi^{0,1}:\Omega^1\to \Omega^{0,1}$ and $\pi^{1,0}:\Omega^1\to \Omega^{1,0}$. Then for
$\xi,\kappa\in \Omega^1$ we define
\[
\<\xi,\overline{\eta}\> = \<\pi^{1,0}\xi,\overline{\pi^{1,0}\eta}\> + \<\pi^{0,1}\xi,\overline{\pi^{0,1}\eta}\> 
\]
or equivalently for the round inner products
\[
(\xi,\eta)=  (\pi^{1,0}\xi,\pi^{0,1}\eta) + (\pi^{0,1}\xi,\pi^{1,0}\eta).
\]

\begin{figure}
\[ \includegraphics[scale=0.85]{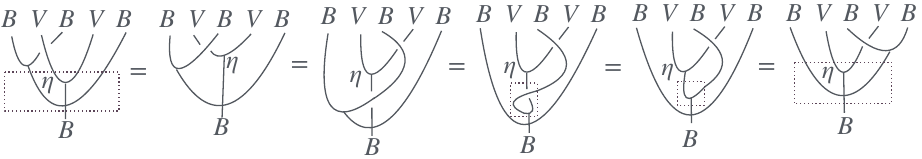}\]
\caption{Proof of  Proposition~\ref{propchernB}. Here $B$ stands for $B_+(V)$. \label{fig3}}
\end{figure}

\begin{proposition}\label{propphim} In the case of $A=B_+(V)$ with $V$ unitary, let $\eta:V\tens V\to B_+(V)$ be any linear map obeying $*\,\eta=\eta\,\dagger$. Then $\phi_-(a\tens v\tens w\tens b)=a\, \eta(v\tens w)b$ and the corresponding $\phi_+$ in (\ref{eq47}) lead to hermitian metrics on $\Omega^{0,1}$ and $\Omega^{1,0}$. A sufficient condition for these maps to descend to $\tens_{B_+(V)}$ is that $\Psi_{V,B_+(V)}$ is involutive and that the image of $\eta$ commutes with $B_+(V)$ via $\Psi$. 
\end{proposition}
\proof  We will write $B:=B_+(V)$ as a shorthand and to align with Section~\ref{secdolb}.  It is immediate that for any linear map $\eta:V\tens V\to B$, 
\begin{equation}\label{phim} \phi_-:=\cdot\circ(\cdot\tens\cdot)(\id\tens\eta\tens \id): \Omega^{0,1}\tens\Omega^{1,0}\to B\end{equation}
is a bimodule map, where the domain is $B\tens V\tens V\tens B$, we apply $\eta$ in the middle and then multiply up in $B$. By the discussion around (\ref{psirev}), we also had
\[ \Theta^{p,q}: \Omega^{0,p}\tens\Omega^{q,0}\to \Omega^{p,0}\tens \Omega^{0,q}\]
as a bimodule isomorphism,  hence 
\[ \phi_+:=\phi_-\circ(\Theta^{1,1})^{-1}:\Omega^{1,0}\tens\Omega^{0,1}\to B\]
is a bimodule map. Similarly, it is immediate that $\phi_-(a\tens v\tens w\tens b)^*=\phi_-(b^*\tens w^*\tens v^*\tens a^*)$ under our assumptions on $\eta$. It is then clear that $\<\, , \>$ on $\Omega^{0,1}$ defined by $\phi_-$ is hermitian. For $\phi_+$ it is similarly clear by the the bimodule map property we need $\phi_+(v\tens 1\tens 1\tens w)=\eta\circ\Psi^{-1}(v\tens w)$ to have the same property as assumed for $\eta$, which follows from (\ref{psidag}). For the last part, the proof that $\phi_-((a\tens v).b\tens w\tens c)=\phi_-(a\tens v)\tens b.(w\tens c))$ is shown in Figure~\ref{fig3} under the further assumption stated. The dotted box is $\phi_-$. The first equality is that $\Psi$ is involutive. The second is that $\eta$ is a morphism so we can move the crossing to the other side. We then use associativity and for the 4th equality we use braided commutativity in the form shown boxed, then associativity again. The  map $\phi_+$ descends to $\tens_B$ if and only if $\phi_-$ does as is clear from (\ref{eq47}) since if either descends, it factors through $\Omega^{1,1}$.  \endproof

Note that while the diagrammatic proof for $\phi_-$ to descend works in the same way if  $\eta:V\tens V\to L$ and hence $\phi_-$ are valued in a bimodule $L$ as in  Section~\ref{secbundle}, this does not help for the quantum plane since $\Psi$ is not involutive.

\subsection{Chern connection for a hermitian metric}\label{secch}

There is a usual definition of hermitian metric preserving connection in noncommutative differential geometry which includes the classical definition when $A$ is a $*$-algebra of complex functions on a manifold. For a left $A$-module $E$, we denote by $\overline{E}$ its conjugate right module (so $\bar e.a=\overline{a^*.e}$ for all $a\in A$ and $e\in E$, and the overline means to view in $\overline{E}$). Then a hermitian inner product is a bimodule map $\<,\>:E\tens\overline{E}\to A$ and a left connection $\nabla_E:E\to \Omega^1\tens_A E$ is hermitian metric compatible if \cite[Def.~8.33]{BegMa}
\[
\extd\, \<e,\overline{f}\> = e_1\,\<e_2,\overline{f}\>+ \<e,\overline{f_2}\> (f_1)^*
\]
for $e,f\in E$, where $\nabla_E(e)=e_1\tens e_2$ is  an explicit notation (sum of such terms understood). Note that this needs $\Omega^1$ to be part of a $*$-calculus. 

For $E$ a left finitely generated projective (f.g.p.) $A$-module we set $E^\flat={}_A\Hom(E,A)$ and take a dual basis $e_i\tens e^i\in E^\flat\tens E$ (summation implicit). This is defined by
$e=\ev(e\tens e_i)\, e^i$ for all $e\in E$. There is a projection matrix with entries in $A$, $P^i{}_j=\ev(e^i\tens e_j)$ associated to the dual basis. If $\<\,,\>:E\tens \overline{E}\to A$ is a hermitian inner product, 
we define $g^{ij}=\<e^i,\overline{e^j}\>$, which is a hermitian matrix with entries in $A$. If the metric is nondegenerate then there is another hermitian matrix with entries in $A$, $\tilde g_{ij}$ so that the matrix product $g^{ij}\,\tilde g_{jk}=P^i{}_k$ and (omitting the matrix indices this time) $\tilde g\, g=P^*$. 

We can also use the dual basis description of the left f.g.p.\  $A$-module $E$ to describe connections on $E$. 
A connection $\nabla_E: E\to \Omega^1\tens_A E$ can be uniquely described by $\nabla_E(e^i)=-\Gamma^i{}_k\tens e^k$ (summation over $k$). (The minus sign is a convention to match with upper and lower index notations for Christoffel symbols in classical geometry.) Here $\Gamma^i{}_k$ is a matrix 
with entries in $\Omega^1{}$. This can be varied for $\partial$-connections and $\bar\partial$-connections by taking entries in $\Omega^{1,0}{}$ and $\Omega^{0,1}{}$ respectively. 

 A left holomorphic connection on $E$ means a left connection of the form $\bar\partial_E:E\to \Omega^{0,1}\tens_A E$ whose holomorphic curvature,
\[
R_{\bar\del_E}=(\bar\partial\tens\id-\id\wedge \bar\partial_E)\bar\partial_E:E\to \Omega^{0,2}\tens_A E
\]
vanishes. One says that $E$ is a holomorphic left $A$-module if it is equipped with such a connection. 
In this context there is a natural construction of metric preserving connections, due to Chern in the classical case.  We present a consequence of Theorems~8.53 and 8.54 in \cite{BegMa}. 

\begin{corollary}\label{cherncor} 
Suppose that $\Omega$ is a factorisable Dolbeault complex, and set $E:=\Omega^{1,0}$. 
\begin{enumerate}
\item $ \bar\partial_E:=\Theta^{0110} \bar\partial:\Omega^{1,0} \to  \Omega^{0,1}\tens_A  \Omega^{1,0}$ is a holomorphic connection making $(E, \bar\partial_E)$ into a holomorphic module. 

\item The  torsion  of $\bar\del_E$  (defined to be $T_{\bar\del_E}=\wedge\bar\del_E-\bar\del: \Omega^{1,0}\to \Omega^{1,1}$) vanishes. 

\item If $E$ is left f.g.p.\ and equipped with a nondegenerate hermitian inner product $\<,\>: E\tens \overline{E}\to A$ then there is a unique left $\partial$-connection $\partial_E:E\to \Omega^{1,0}\tens_A E$ on $E$ such that the ordinary connection
\[ \nabla_E=\bar\partial_E+ \partial_E:\Omega^{1,0}\to \Omega^1\tens_A\Omega^{1,0}\]
is  hermitian metric compatible. If we write the Christoffel symbols of the given $ \bar\partial_E$ as the matrix $\Gamma_-$ then the Christoffel symbols of $ \partial_E$ are
\[
\Gamma_+= - \partial(g)\,\tilde g - g\, (\Gamma_-)^*\,\tilde g
\]
where $(\Gamma_-)^*$ is the star transpose of the matrix. 

\item $ \bar\partial_E$ is a bimodule connection with $\sigma_E{}^{0,1}=-\Theta^{0110}(\wedge):\Omega^{1,0}\tens_A\Omega^{0,1}  \to 
\Omega^{0,1}\tens_A\Omega^{1,0}$.

\item  If  $\< \, ,\>$ descends to
$\<\,,\>:\Omega^{1,0}\tens_A\overline{\Omega^{1,0}}\to A$ then $\partial_E$ and $\nabla_E$ are bimodule connections with  $\sigma_E{}^{1,0}:\Omega^{1,0}\tens_A\Omega^{1,0}  \to 
\Omega^{1,0}\tens_A\Omega^{1,0} $ and $\sigma_E=\sigma_E{}^{0,1}+\sigma_E{}^{1,0}$ respectively.
\end{enumerate}
\end{corollary}
\proof First we check that $\bar\partial_E:\Omega^{1,0} \to  \Omega^{0,1}\tens_A  \Omega^{1,0}$ is a $\bar\partial$-connection, 
\[
 \bar\partial_E(a\,\xi)=\Theta^{0110} (\bar\partial a\wedge\xi+ a\, \bar\partial \xi)= \bar\partial a\tens\xi + a\,  \bar\partial_E(\xi). 
\]
Its curvature is $R_{\bar\del_E}(\xi)= \bar\partial(\kappa)\tens\eta-\kappa\wedge \bar\partial_E(\eta)$, where $\bar\partial_E(\eta)=\kappa\tens\eta$. It follows that $\wedge R_{\bar\del_E}(\xi)=  \bar\partial^2\xi=0\in \Omega^{1,2}$ and so $R_{\bar\del_E}(\xi)=\Theta^{0210}(\bar\partial^2\xi)=0$. Hence, we have a holomorphic bundle. Now we apply the results in \cite{BegMa} to construct $\partial_E$ as stated. Note that  applying $*$ to the entries of $\Gamma_-$ lands them in $\Omega^{1,0}$. Moreover, as $\Theta^{0110} $ is inverse to $\wedge$, we have  $ \bar\partial_E$ torsion free. \endproof
 
 Note that the curvature of $\nabla_E$ restricts to $R_E:\Omega^{1,0}\to \Omega^{1,1}\tens_A\Omega^{1,0}$ and, as the torsion of $\bar\partial_E$ is zero, the torsion of $\nabla_E$ restricts to $T_E:\Omega^{1,0}\to \Omega^{2,0}$. 

\begin{remark} In the Chern construction, the left $\bar\partial$-connection $ \bar\partial_E:
E\to   \Omega^{0,1}\tens_A E$ is equivalent to a right $\partial$-connection $\del_{\overline{E}}^R: \overline{E}\to \overline{E}\tens_A\Omega^{1,0}$
by $\del_{\overline{E}}^R(\overline{e})=\overline{f}\tens\xi^*$ where $\bar\partial_E e=\xi\tens f$ (summation implicit). 
 Also note that we do not assume that $\<\, , \>$ descends to $\tens_A$, though it may of course do so. 
  \end{remark}

We give two classes of examples, the first is the general quantum-braided plane. 

\begin{proposition}\label{propchernB} In the setting of Proposition~\ref{propphim} for $A=B_+(V)$ with $V$ a unitary object, we take $v^i$ to be a basis of $V$ so that $e^i=v^i\tens 1$ is a basis over the algebra of $E=\Omega^{1,0}$. The hermitian metric on $E=\Omega^{1,0}$ has 
\[
g_E^{ij}=\<e^i,\overline{e^j}\>_E= \phi_+(e^i\tens e^{j*})=\phi_+(v^i\tens 1\tens 1 \tens v^j{}^*)= 
\eta(\Psi^{-1}(v^i\tens v^j{}^*)),
\]
where we require $\eta(v^i\tens v^{j*})=\eta(v^j\tens v^{i*})^*$. If the matrix $g_E$ is invertible with inverse
$\tilde g_E$ then the Chern connection on $E=\Omega^{1,0}$
is given by
\[
\bar\del_E e^i=0,\quad \del_E e^i = \del(g_E^{ij})\, \tilde g_{Ejk} \tens e^k
\]
and the torsion of $\nabla_E=\bar\del_E+\del_E$ by
\[
T_E (e^i )= \del(g_E^{ij})\, \tilde g_{Ejk} \wedge e^k.
\]
 \end{proposition}
\proof  By the formulae in Proposition~\ref{propdoublecomp} both $\del$ and $\bar\del$ applied to $e^i$ give zero as the $B_+(V)$-component is $1$. 
As $\bar\del e^i=0$ we get $\bar\del_E e^i=0$ as the first part of the Chern connection. Then Corollary~\ref{cherncor} gives
the formula for $ \del_E e^i $
and the formula for the torsion is immediate.
\endproof
Concrete examples on the quantum plane $\C_q^2$ will be given later. Our second class of examples 
are based on a discrete group $G$ with left-invariant calculus corresponding to a Cayley graph generated by a subset $\CC\subseteq G\setminus\{e\}$. The algebra $A=C(G)$ has free left invariant generators $e^a$ for $a\in \CC$ and 
\[
e^a\, f=R_a(f)\, e^a\,\quad \extd f=\sum_{a\in \CC} (R_a(f)-f)\, e^a,
\]
where $R_a(f)(g)=f(ga)$ for $f\in C(G)$. 
If $\CC$ is closed under inverse then this is a star calculus with $e^a{}^*=-e^{a^{-1}}$, and it is bicovariant if $\CC$ is a union of conjugacy classes. 
There are natural exterior algebras 
\[ \Omega_L(G) \twoheadrightarrow \Omega_{LL}(G)\twoheadrightarrow \Omega_{\mathrm{Wor}}(G)\]
defined as successive quotients of the maximal prolongation by relations \cite[Cor.~1.54]{BegMa}
\begin{equation}\label{eqnLL} \sum_{a,b\in \CC; ab=z} e^a\wedge e^b=0\end{equation}
for all $z\in G\setminus\{e\}$ for $\Omega_L$ and including $z=e$ for $\Omega_{LL}$. The well-known case here is the  Woronowicz calculus $\Omega_{\mathrm{Wor}}$ which applies in the bicovariant case and  has relations in degree 2 given by setting to zero elements invariant under the braiding
\[ \Psi(e^a\tens e^b)= e^{aba^{-1}}\tens e^a\]
for all $a,b\in \CC$. This is immediate from the general construction for any Hopf algebra in \cite{Wor}. 

\begin{proposition} \label{propGdolb}
 If $\CC=\CC^{1,0}\sqcup C^{0,1}$ is a disjoint union of two subsets bijective by group inversion then $\Omega_L(G),\Omega_{LL}(G)$ and in the bicovariant case $\Omega_{\rm Wor}(G)$ have quotients $\Omega(G)$ which form a  Dolbeault complex. Moreover, the further relations
\begin{align} \label{hur}
e^a\wedge e^{b^{-1}}+\wedge\Psi(e^a\tens  e^{b^{-1}}) =0,
\end{align}
for all $a,b\in \CC^{0,1}$ (if not already present) make $\Omega(G)$  factorisable in degree 2. 
 \end{proposition}
\proof  (1) The quotient procedure for $\Omega(G)$ is to take each of the relations of $\Omega_{L}(G)$, $\Omega_{LL}(G)$ and $\Omega_{\rm Wor}(G)$ and project to each $\Lambda^{p,q}$ component as separate relations.  Next, the calculus  $\Omega_L(G)$ is known \cite[Prop.~1.53]{BegMa} to be inner to all orders by the element $\theta=\sum_{b\in\CC} e^b$ so that that $\extd$ is the graded commutator 
\[ 
\extd \xi =[\theta,\xi \} = \theta\wedge\xi-(-1)^{|\xi|}\xi\wedge\theta. 
\]
In our case, we split $\theta=\theta^{0,1}+\theta^{1,0}$ where $\theta^{0,1}=\sum_{a\in \CC^{0,1}}e^a$ and correspondingly for $\theta^{1,0}$. Then $\del = [\theta^{1,0},\  \}$ and $\bar\del = [\theta^{0,1},\  \}$ sum to $\extd$, and $\del:\Omega^{p,q}\to \Omega^{p+1,q}$ and $\bar\del:\Omega^{p,q}\to \Omega^{p,q+1}$, showing the integrability of the complex structure. Explicitly, on degree 0, $\bar\del f= \sum_{a\in \CC^{0,1}}\del_a (f) e^a$ and $\del f= \sum_{a\in \CC^{1,0}}\del_a (f) e^a$ are the two pieces of $\extd f$ above.  The $*$-structure of $\Omega_L(G)$  is given by $(e^a)^*=-e^{a^{-1}}$ which takes us to the other side. These constructions respect the total grading and hence descend to our quotient $\Omega(G)$. Similarly  for the two other cases. This constructs the Dolbeault complex based on this  $\Omega(G)$.  

(2) $\Omega(G)$ thus far need not factorisable. We show that we can add the additional relations as stated to make it so at least in degree 2. It is enough to show this when starting with $\Omega_{\rm Wor}(G)$ as $\Omega_L(G)$ and $\Omega_{LL}(G)$ have finer relations. 
In degree 2
$R_{\mathrm{Wor}} =\ker(\id-\Psi)\subset V^{2,0}\oplus V^{1,1}\oplus V^{0,2}$. Since each of these direct summands is invariant to $\Psi$, taking the kernel and then projecting to $V^{1,1}$ is the same as taking $\ker(\id-\Psi)$ for the map on $V^{1,1}$, which we do. Now set $V^{1001}=
\Lambda^{1,0}\tens\Lambda^{0,1} $ and $V^{0110}=
\Lambda^{0,1}\tens\Lambda^{1,0} $. We describe elements of $V^{1,1}$ as 
$\left(\begin{smallmatrix}x\\y\end{smallmatrix}\right)$ where $x\in V^{1001}$ and $y\in V^{0110}$ and write $\Psi=
\left(\begin{smallmatrix}0&b\\c&0\end{smallmatrix}\right)$. Then, restricting to $V^{1,1}$, we get the Woronowicz relations.
\[
\ker(\id-\Psi) = \big\{ \left(\begin{smallmatrix}z\\ cz \end{smallmatrix}\right) \,:\, z\in  V^{1001}\ \mathrm{and}\  bcz=z\big\}.
\]
We add relations of the form
\[
(\id-\Psi) V^{1001}= \big\{ \left(\begin{smallmatrix}x\\ cx \end{smallmatrix}\right) \,:\, x\in  V^{1001} \big\},
\]
and note that the additional relations contain the Woronowicz relations, so we just work with them. We then have:
\newline i)\quad The intersection of $(\id-\Psi) V^{1001}$ with either $V^{1001}$ or $V^{0110}$ is zero; we need this to ensure that the 
$V^{1001}$ or $V^{0110}$ are not quotiented by the relations. 
\newline ii)\quad Given the relations, we have $\left(\begin{smallmatrix}x\\ 0 \end{smallmatrix}\right)\cong -\left(\begin{smallmatrix}0\\ cx \end{smallmatrix}\right)$ with $c$ invertible, completing factorisability in degree 2.
\endproof

We are also free to make further quotients of each holomorphic and holomorphic factor as long as we preserve the Dolbeault complex structure.  

\begin{proposition}\label{propchernG} Let $\Omega(G)$ be a factorisable Dolbeault complex of the form in the preceding proposition and lemma with both $\CC^{1,0},\CC^{0,1}$ conjugacy classes. We suppose that we are given a nondegenerate hermitian metric of the form $\<e^{a^{-1}}, e^{b^{-1}}\>=g^{a^{-1}b^{-1}}\in C(G)$ on $E=\Omega^{1,0}$, where $a,b\in \CC^{0,1}$. Here $E$ is spanned over $A=C(G)$ by $\{e^{a^{-1}}\}$ for $a\in \CC^{0,1}$. Then we have a metric preserving connection on this given by
\begin{align*}
\nabla_E e^{b^{-1}} &= \sum_{a\in  \CC^{0,1}} e^a\tens \big(e^{b^{-1}} - e^{a^{-1}b^{-1}a} \big)  \cr
&\quad 
 + \sum_{a\in\CC^{0,1}}   e^{a^{-1}} \tens e^{b^{-1}}
  -   \sum_{a,s,c\in\CC^{0,1}}   g^{b^{-1} , a^{-1}s^{-1}a}  \, e^{a^{-1}}     \,\tilde g_{   s^{-1}c^{-1}}  \tens e^{c^{-1}}
\end{align*}
 and its torsion is
  \[
T_{E} (e^{b^{-1}}) =-\sum_{a\in\CC^{0,1}}    e^{b^{-1}} \wedge e^{a^{-1}} 
  -   \sum_{a,s,c\in\CC^{0,1}}   g^{b^{-1} , a^{-1}s^{-1}a}  \, e^{a^{-1}}     \,\tilde g_{   s^{-1}c^{-1}}  \wedge e^{c^{-1}}
 \]
(Note that in this case $\tilde g_{bc}$ is the matrix inverse of $g^{ab}$.) Moreover, the hermitian metric descends to $\tens_A$ if and only if 
the matrix $g^{ab}$ is diagonal, and in this case we have a bimodule connection, with for $a\in \CC^{0,1}$,
 \begin{align*}
\sigma_E(e^{b^{-1}} \tens e^a)&=e^a\tens e^{a^{-1}b^{-1}a},\\
  \sigma_E(e^{b^{-1}}\tens e^{a^{-1}})&= g^{b^{-1}, b^{-1}}    \, R_{a^{-1}}(\tilde g_{ab^{-1}a^{-1}, ab^{-1}a^{-1}  }) e^{a^{-1}} \tens e^{ab^{-1}a^{-1} }.
\end{align*}
\end{proposition}
\proof Begin with
\[
\bar\del e^{b^{-1}}=\sum_{a\in  \CC^{0,1}} \big(e^a\wedge e^{b^{-1}} + e^{b^{-1}} \wedge e^a\big) = 
\sum_{a\in  \CC^{0,1}} e^a\wedge \big(e^{b^{-1}} - e^{a^{-1}b^{-1}a} \big) 
\]
for $b^{-1}\in  \CC^{1,0}$, where we have used the relation (\ref{hur}). Thus we have the first term (the $\bar\del$-connection $\bar\del_E$) in the equation for $\nabla_E e^{b^{-1}}$
in the statement. We check that this is a holomorphic connection (which would be automatic if we had factorisation to degree 3). The holomorphic curvature is
\begin{align*}
R_{\mathrm{Hol}}(e^{b^{-1}})&=\sum_{a\in  \CC^{0,1}} \bar\del e^a\tens \big(e^{b^{-1}} - e^{a^{-1}b^{-1}a} \big) \cr
&\quad -\sum_{a,c\in  \CC^{0,1}} e^a\wedge e^c\tens \big(e^{b^{-1}} - e^{a^{-1}b^{-1}a} - e^{c^{-1}b^{-1}c} 
+  e^{c^{-1}a^{-1}b^{-1}ac}\big) \cr
&= \sum_{a,c\in  \CC^{0,1}} e^a\wedge e^c\tens \big(e^{b^{-1}} 
-  e^{(ac)^{-1}b^{-1}ac}\big)=0
\end{align*}
 by looking at the sums where $ac$ is a constant and using the definition of $\Omega_L$ (the weakest relations we consider). 

Then $\bar\del_E$ is a bimodule $\bar\del$-connection and an explicit calculation using
\begin{align}  \label{pfj}
\sigma_E( e^{b^{-1}} \tens \bar\del f)= \bar\del_E (e^{b^{-1}} \,f) -  \bar\del_E (e^{b^{-1}}) \,f
\end{align}
 gives the first part of $\sigma_E$ as stated.  
The Christoffel symbols for $\bar\del_E$ are
\[
(\Gamma_-)^{b^{-1}}{}_{c^{-1}} = - \delta_{b,c} \sum_{a\in\CC^{0,1}} e^a  + \sum_{a\in\CC^{0,1}: c^{-1}=a^{-1}b^{-1}a} e^a  
\]
and then the Christoffel symbols of $ \partial_E$ are
\begin{align*}
- (\Gamma_+)^{b^{-1}} {}_{c^{-1}}&= \big( \partial(g^{b^{-1} s^{-1}} )    + g^{b^{-1} u^{-1}} (\Gamma_-)^*{}_{u^{-1}}{}^{s^{-1}}
\big)    \,\tilde g_{   s^{-1}c^{-1}} \cr
&= \big( \partial(g^{b^{-1} s^{-1}} )  + g^{b^{-1} u^{-1}}\, \delta_{u,s} \sum_{a\in\CC^{0,1}} e^{a^{-1}}
  - g^{b^{-1} u^{-1}}  \sum_{a\in\CC^{0,1}: u^{-1}=a^{-1}s^{-1}a} e^{a^{-1}}  
\big)    \,\tilde g_{   s^{-1}c^{-1}} \cr
&= \big( \sum_{a\in\CC^{0,1}} e^{a^{-1}}\,  g^{b^{-1} s^{-1}}
  - \sum_{a\in\CC^{0,1}: u^{-1}=a^{-1}s^{-1}a}     g^{b^{-1} u^{-1}}  \, e^{a^{-1}}  
\big)    \,\tilde g_{   s^{-1}c^{-1}} \cr
&=\delta_{b,c}\sum_{a\in\CC^{0,1}}   e^{a^{-1}}
  -   \sum_{a\in\CC^{0,1}}   g^{b^{-1} , a^{-1}s^{-1}a}  \, e^{a^{-1}}     \,\tilde g_{   s^{-1}c^{-1}}, 
\end{align*}
where we have used that $ \sum_{a\in\CC^{0,1}} e^{a^{-1}}$ is the inner element defining $\bar\del$. This gives the second line for $\del_E$
in the formula for $\nabla_E$  in the statement. Also as the algebra is commutative, if the metric descends to $\tens_A$ we can only have nonzero values of 
 $\<e^a,\overline{e^b}\>$ if $ab^{-1}=e$, so in this case $g^{ab}$ is diagonal and we get
\begin{align*}
\del_E e^{b^{-1}} &=  \sum_{a\in\CC^{0,1}}   e^{a^{-1}} \tens e^{b^{-1}}
  -   \sum_{a\in\CC^{0,1}}   g^{b^{-1} , b^{-1}}    \, R_{a^{-1}}(\tilde g_{   ab^{-1}a^{-1}   ,   ab^{-1}a^{-1}  })    \, e^{a^{-1}}     \tens e^{ab^{-1}a^{-1} }.
\end{align*}
 Now using (\ref{pfj}) for $\del$ instead of $\bar\del$ gives the second part of $\sigma_E$. \endproof

\begin{corollary} If $G$ is Abelian then the above calculus with the Woronowicz calculus on each factor reduces to $\Omega(G)=\Omega_{\mathrm{Wor}}(G)$ as the Grassmann algebra on the  basic 1-forms $e^a,e^{b^{-1}}$. Moreover, if  $g^{a^{-1}b^{-1}}=\delta_{a,b}$ is the constant  `Euclidean metric' then the hermitian metric compatible connection $\nabla_E$ in Proposition~\ref{propchernG} is torsion free.
\end{corollary} 
\proof The first part is immediate as the conjugations are absent and $\Psi$ is then the flip map. For the second part, for the constant metric the torsion then becomes $-\sum_{a\in C^{0,1}}(e^{b^{-1}}\wedge e^{a^{-1}}+e^{a^{-1}}\wedge^{b^{-1}})=0$. \endproof

We similarly have the Grassmann algebra and zero torsion for any group and any diagonal metric if $C^{0,1}$ is a singleton set. 

\subsubsection{Factorisable complex structure on alternating group $A_4$}

The normal subgroup $G=A_4\subset S_4$ of even permutations on 4 elements is an order 12 solvable group which we take with generators 
\[ t=(123),\ u=(14)(23),\ v=(12)(34),\ w=(13)(24),\]
\[ x=tv=ut=(134),\ y=tw=vt=(243),\ z=tu=wt=(142).\]
Here $t,x,y,z$ are all of order 3 and $u,v,w$ are of order 2. It as already observed in \cite{BegMa} that $\CC^{0,1}=\{t,x,y,z\}$ is a conjugacy class disjoint from its inverses $\CC^{1,0}=\{t^{-1},x^{-1},y^{-1},z^{-1}\}$ and hence $\CC=\CC^{0,1}\amalg\CC^{0,1}$ and Proposition~\ref{propGdolb}  allows for the construction of a factorisable double complex as follows. One can check that $\CC^{1,0}\CC^{1,0}=\CC^{0,1}$, $\CC^{0,1}\CC^{0,1}=\CC^{1,0}$ and $\CC^{0,1}\CC^{1,0}=\CC^{1,0}\CC^{0,1}=\{e,u,v,w\}$ (so in fact all elements of the group are obtained in two steps from $e$ on the Cayley graph. 

Starting with $\Omega_L(G)$, there are $8$ products of the form $ab=z$ for $a,b\in\CC$ if $z=e,u,v,w$, so that (\ref{eqnLL}) for the latter three gives respectively the 8-term relations (omitting $\wedge$ to save space)
\[ e^x e^{t^{-1}}+ e^{y^{-1}}e^x+e^z e^{y^{-1}}+ e^{t^{-1}}e^z+ e^t e^{x^{-1}}+ e^{z^{-1}}e^t+ e^y e^{z^{-1}}+ e^{x^{-1}}e^y=0\]
\[  e^t e^{y^{-1}}+e^{x^{-1}}e^t+ e^z e^{x^{-1}}+e^{y^{-1}}e^z+  e^x e^{z^{-1}}+e^{t^{-1}}e^x+ e^y e^{t^{-1}}+e^{z^{-1}}e^y=0\]
\[ e^t  e^{z^{-1}} + e^{y^{-1}}e^t+ e^x e^{y-1}+e^{z^{-1}}e^x+ e^y e^{x^{-1}}+e^{t^{-1}}e^y+ e^z e^{t^{-1}}+e^{x^{-1}}e^z=0\]
for $\Omega_L(G)$ and we have additonally
the 8 term relation
\[ e^te^{t^{-1}}+e^{t^{-1}}e^t+  e^xe^{x^{-1}}+e^{x^{-1}}e^x+e^ye^{y^{-1}}+e^{y^{-1}}e^y+e^ze^{z^{-1}}+e^{z^{-1}}e^z=0\]
for $\Omega_{LL}(G)$. In addition, for either case, we have the eight 4-term relations for $ab=z$, $z=t,x,y,z$ and their inverses. These are in order 
\[ e^{t^{-1}}e^{t^{-1}}+e^{x^{-1}} e^{z^{-1}}+e^{z^{-1}} e^{y^{-1}}+ e^{y^{-1}} e^{x^{-1}}=0,\quad e^{x^{-1}}e^{x^{-1}}+e^{t^{-1}} e^{y^{-1}}+e^{y^{-1}} e^{z^{-1}}+e^{z^{-1}} e^{t^{-1}}=0 \]
\[ e^{z^{-1}}e^{z^{-1}}+e^{t^{-1}} e^{x^{-1}}+ e^{x^{-1}} e^{y^{-1}}+ e^{y^{-1}} e^{t^{-1}}=0,\quad e^{y^{-1}}e^{y^{-1}}+ e^{t^{-1}} e^{z^{-1}}+ e^{z^{-1}} e^{x^{-1}}+ e^{x^{-1}} e^{t^{-1}}=0,\quad \]
\[e^t e^t+e^x e^y+e^y e^z+ e^z e^x=0,\quad e^x e^x +e^t e^z+e^z e^y+e^y e^t =0,\]
\[ e^y e^y+e^t e^x+ e^x e^z+ e^z e^t=0,\quad e^z e^z+e^t e^y+e^y e^x+e^x e^t=0.\]

By contrast, $\Omega_{\mathrm{Wor}}(G)$ in degree 2 has each half vanishing of each of the three 8-term relation of $\Omega_{L}(G)$ and it has 12  further relations
 \[ (e^a)^2=0, \quad \{e^a,e^{a^{-1}}\}=0,\quad  \forall a\in \CC\]
 These render the additional relation of $\Omega_{LL}(G)$ redundant and the eight 4-terms relations of $\Omega_{LL}(G)$ into eight 3-term relations, These square and 3-term relations are exactly $\Omega_{\rm Wor}^{0,\bullet}$ and $\Omega_{\rm Wor}^{\bullet,0}$ in degree 2 taken separately (note that these are known to have further relations in higher degree\cite{NML}). This summarises the initial exterior algebras for our construction. Here $\Omega^2_L(G)$ is 53-dimensional, $\Omega^2_{LL}(G)$ is 52-dimensional and $\Omega^2_{\rm Wor}(G)$ is 38-dimensional.
 
 Now applying the first part of Proposition~\ref{propGdolb}, all of the relations (in degree 2) are contained entirely in one of spaces $\Omega^{2,0},\Omega^{1,1},\Omega^{0,2}$ so there is no need for any division of the relations, hence all $\Omega_L(G),\Omega_{LL}(G), \Omega_{\rm Wor}(G)$ are already the desired Dolbeault complex $\Omega(G)$. Graded anticommutator with $\bar\theta=e^t+e^x+e^y+e^z$ gives $\bar\del$ and with $\theta=e^{t^{-1}}+e^{x^{-1}}+e^{y^{-1}}+ e^{z^{-1}}$ gives $\del$ in the double complex. This is not in general, however, factorisable as there will not be enough relations in $\Omega^{1,1}$ to put all elements in a given order or say holomorphic wedge antiholomorphic. For this, we impose the further relations displayed in Proposition~\ref{propGdolb}. This splits each of the 8-term relations into 4 pairs. Hence the factorisable $\Omega^2(G)$  is 44-dimensions coming from $\Omega_L(G)$, $41$-dimensional coming from $\Omega_{LL}(G)$ and 32-dimensional coming from $\Omega_{\rm Wor}(G)$. 
 
One can then apply the Chern construction here (and later, the doubled  Chern construction).  The torsion, however, typically will not vanish. For example, the  Euclidean metric $g^{ab}=\delta_{a,b}$ for all $a,b\in\CC$ gives
 \[ T_E(e^{a^{-1}})= \theta e^{a^{-1}} \]
for all $a\in \CC^{1,0}$  if we start with $\Omega_{\rm Wor}(G)$. A framed quantum geometry on this $G$ with $\Omega_{\rm Wor}(G)$ (using a different formalism) is in \cite{NML}.

\subsection{Doubled Chern construction for hermitian metrics on $\Omega^1$}\label{secdoublech}
In Corollary~\ref{cherncor}, we considered the Chern construction from the usual point of view, with holomorphic bundles having a 
$\bar\partial$-connection. However, we could swap the roles of $\Omega^{0,1}$ and $\Omega^{1,0}$ in the theory, which would be equivalent to writing usual complex analysis in terms of functions of $\bar z\in\C$ rather than $z\in\C$. If we do this then we get a conjugate version of Corollary~\ref{cherncor}, which we can summarise as follows:

\begin{corollary}\label{anti-cherncor} Let $\Omega$ be a factorisable Dolbeault complex and set $F:=\Omega^{0,1}$. 

(1) $ \partial_F:=\Theta^{1001} \partial:\Omega^{0,1} \to  \Omega^{1,0}\tens_A  \Omega^{0,1}$ is a torsion free anti-holomorphic connection making $(F, \partial_F)$ into a anti-holomorphic module. 

(2) If $F$ is left f.g.p.\ and equipped with a nondegenerate hermitian inner product $\<,\>: F\tens \overline{F}\to A$, there is a unique left $\bar\partial$-connection $\bar\partial_F:F\to \Omega^{0,1}\tens_A F$ on $E$ such that the ordinary connection
\[ \nabla_F=\bar\partial_F+ \partial_F:\Omega^{0,1}\to \Omega^1\tens_A\Omega^{0,1}\]
is  hermitian metric compatible. If we write the Christoffel symbols of $ \partial_F$ as the matrix $\Gamma_+$ then the Christoffel symbols of $ \bar\partial_F$ are
\[
\Gamma_-= - \bar\partial(g)\,\tilde g - g\, (\Gamma_+)^*\,\tilde g.
\]

(3) In this case,  $ \bar\partial_F$ is a bimodule connection with 
\[\sigma_F{}^{1,0}=-\Theta^{1001}(\wedge):\Omega^{0,1}\tens_{\blu A}\Omega^{1,0}  \to 
\Omega^{1,0}\tens_{\blu A}\Omega^{0,1}.\]

(4) Moreover,  if  $\< \, ,\>$ descends to
$\<\,,\>:\Omega^{0,1}\tens_A\overline{\Omega^{0,1}}\to A$ then $\bar\partial_F$ and $\nabla_F$ are bimodule connections with respective generalised braidings 
\[ \sigma_F{}^{0,1}:\Omega^{0,1}\tens_{\blu A}\Omega^{0,1}  \to 
\Omega^{0,1}\tens_{\blu A}\Omega^{0,1},\quad \sigma_F=\sigma_F{}^{0,1}+\sigma_F{}^{1,0}.\]
\end{corollary}

In particular, we have the following conjugate version of Proposition~\ref{propchernB}.

\begin{lemma}\label{propchernC} In the setting of Proposition~\ref{propphim} for $A=B_+(V)$ with $V$ a unitary object and basis $v^i$ to be a basis of $V$, we let $f^i=1\tens v^i$ be a basis of $F=\Omega^{0,1}$.  The hermitian metric on $F=\Omega^{0,1}$ has
\[
g_F^{ij}=\<f^i,\overline{f^j}\>_F = \phi_-(f^i\tens f^{j*})=\phi_-(1\tens v^i\tens v^{j*}\tens 1)=\eta(v^i\tens v^{j*}),
\]
where  $\eta(v^i\tens v^{j*})=\eta(v^j\tens v^{i*})^*$ as before. If the matrix $g_F$ is invertible with inverse
$\tilde g_F$ then the Chern connection on $F=\Omega^{0,1}$
is given by
\[
\del_F f^i=0,\quad \bar\del_F f^i = \bar\del(g_F^{ij})\, \tilde g_{Fjk} \tens f^k
\]
and the torsion of $\nabla_F=\bar\del_F+\del_F$ is 
\[
T_F (f^i )=  \bar\del(g_F^{ij})\, \tilde g_{Fjk} \wedge f^k.
\]
 \end{lemma}
\proof  By the formulae in Proposition~\ref{propdoublecomp} both $\del$ and $\bar\del$ applied to $f^i$ give zero as the $B$-component is $1$. 
As $\del f^i=0$ we get $\del_F f^i=0$ as the first part of the Chern connection. Then Corollary~\ref{cherncor} (reversed) gives 
the formula for $\bar\del_F f^i $
and the formula for the torsion is immediate.
\endproof

Now suppose that there is a Hermitian inner product on $\Omega^1=\Omega^{0,1}\oplus\Omega^{1,0}$ for which
$\Omega^{0,1}$ and $\Omega^{1,0}$ are perpendicular. Then we can use Corollary~\ref{cherncor} to obtain a connection
$\nabla_E:E\to \Omega^1\tens_A E$ for $E=\Omega^{1,0}$ and Corollary~\ref{anti-cherncor} to obtain a connection
$\nabla_F:F\to \Omega^1\tens_A F$ for $F=\Omega^{0,1}$. The direct sum $\nabla_E\oplus\nabla_F$ is then a metric preserving connection on $E\oplus F=\Omega^1$. In particular, we have the following:

\begin{corollary} \label{doubleB}
For $A=B_+(V)$ with $V$ a unitary object, the hermitian metrics 
on $E=\Omega^{1,0}$ in Proposition~\ref{propchernB} and $F=\Omega^{0,1}$ in 
Lemma~\ref{propchernC} given by
$\eta:V\tens V\to B_+(V)$ extend to a hermitian metric on $\Omega^1=E\oplus F$ by taking $E$ and $F$ to be perpendicular. 
 Then the left connection $\nabla=\nabla_{E\oplus F}$ on $\Omega^1$ restricting to $\nabla_E$ and $\nabla_F$ on the two parts is hermitian metric compatible for $\<\, , \>$ and has torsion restricting to $T_E$ and $T_F$. Moreover, $\nabla$ is a bimodule connection if $\<\, ,\>$ descends to $\tens_{B_+(V)}$.
\end{corollary}
\proof Hermitian inner products and hermitian metric preservation were recalled in 
Section~\ref{secch}. From these, the properties on $\Omega^1=E\oplus F$ are proved by considering all choices of restricting the inputs to $E$ or $F$. For example, checking metric preservation 
for $\xi\in E$ and $\kappa\in F$ results in
\[
\xi_1\, \<\xi_2,\overline{\kappa}\> + \<\xi,\overline{\kappa_1}\>\,\kappa_2{}^*=\extd\<\xi,\overline{\kappa}\>
\]
where $\nabla\xi=\nabla_E\xi=\xi_1\tens\xi_2\in \Omega^1\tens_A E$ and $\nabla\kappa=\nabla_F\kappa=\kappa_1\tens\kappa_2\in \Omega^1\tens_A F$ (sum implicit), and this is true as all the inner products in the equation are zero,
because $E$ and $F$ are perpendicular. 
 \endproof

  \begin{proposition} \label{propBdouble}
(1) For $A=B_+(V)$, the unique connection $\nabla_F:F\to\Omega^1\tens_A F$ on $F=\Omega^{0,1}$ that vanishes on all $\bar\del v$ for $v\in V$ is a torsion free bimodule connection with 
\[ \sigma_F(\bar\del v\tens \extd u)=(\extd \tens \bar\del )\,\Psi(v\tens u)\]
and is the Chern connection for any $\eta$ taking values in the complex numbers.
\newline (2)  The same holds for the unique connection  $\nabla_E:E\to\Omega^1\tens_A E$ on $E=\Omega^{1,0}$ that vanishes on all $\del v$ for $v\in V$, now with 
\[ \sigma_E(\del v\tens \extd u)=(\extd\tens \del )\,\Psi^{-1}(v\tens u).\]
\newline (3) If both (1) and (2) apply then 
\begin{align} \label{tydo}
 \sigma_E\,\dag\,\sigma_F=\dag,\quad \nabla_E=\sigma_E\,\dag\,\nabla_F\,*
\end{align}
and   $\nabla=\nabla_E\oplus\nabla_F$ on $\Omega^1=E\oplus F$ is a $*$-preserving bimodule connection.
\end{proposition}
\begin{figure}
\[ \includegraphics{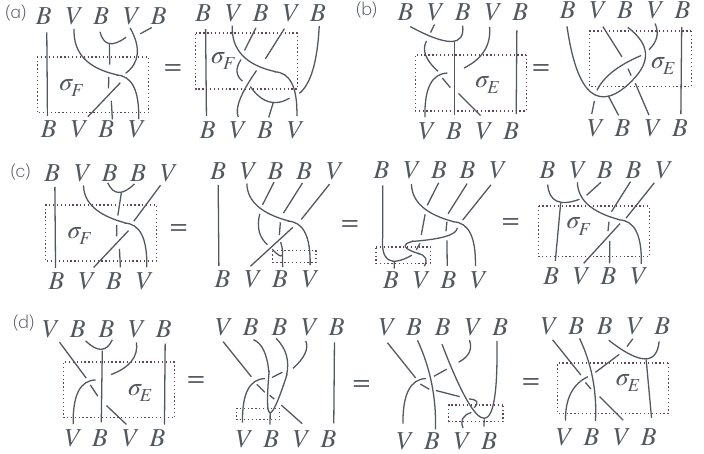}\]
\caption{Diagrams for the proof of Proposition~\ref{propBdouble}. Here $B$ stands for $B_+(V)$. \label{figsigma}}
\end{figure}
\proof  (1) We write $B:=B_+(V)$ for brevity when referring to the braided constructions. Here,  $\Omega^{0,1}=B\tens V$ and $\Omega^{1,0}=V\tens B$ in the category with
$\bar\del v=1\tens v$ and $\del v=v\tens 1$. Using the braided tensor product for the module structure, we can write the relations on the 1-forms
as
\begin{align} \label{pruo}
\bar\del v\, u =\cdot \,(\id\tens\bar\del)\,\Psi\,(v\tens u)\,\quad 
v\, \del u =\cdot \,(\del\tens\id)\,\Psi\,(v\tens u)
\end{align}
 for $v,u\in V$. Then, by definition
\[
\sigma_F(\bar\del v\tens\extd u)=\nabla_F(\bar\del v\, u)-\nabla_F(\bar\del v)\, u = 
\nabla_F(\cdot \,(\id\tens\bar\del)\,\Psi\,(v\tens u))=(\extd \tens \bar\del )\,\Psi(v\tens u)
\]
provided this extends to a bimodule map $\sigma_F:\Omega^{0,1} \tens_B \Omega^{1}  \to \Omega^{1} \tens_B \Omega^{0,1}$. Part of this is the map $\Theta^{1,1}$ in Corollary~\ref{corfac} which necessarily descends to $\tens_B$ (at least in the factorisable case as here). The remaining part part of the required $\sigma_F$ is the map $\Omega^{0,1} \tens_B \Omega^{0,1}  \to \Omega^{0,1} \tens_B \Omega^{0,1}$ shown in the dashed box in Figure~\ref{figsigma}, where part (a) checks that it is a right module map and (c) checks that it is well-defined over $\tens_B$.  The middle step
is to view the small dashed box as the left action of $B$ on the second copy of $\Omega^{0,1}$ and reverse it to the right action the other copy.  Similarly for the corresponding result for $E$ we need to rewrite the last equation in (\ref{pruo}) as 
\begin{align*} 
\del v\, u =\cdot \,(\id\tens\del)\,\Psi^{-1}\,(v\tens u).
\end{align*}
and Figure~\ref{figsigma} shows how the inverse braiding on $V$ here extends via $\sigma_E$ defined as $(\Theta^{1,1})^{-1}$ and the bimodule map shown in the dashed box. Part (b) checks that it is a left module map and part (d) that it descends to the $\tens_B$. The middle step
is to view the small dashed box as the right action of $B$ on the first copy of $\Omega^{1,0}$ and reverse it to the left action the other copy. 

Next, by definition $\extd(\bar\del v)=0$, so by factorisability we have $\del_F:F\to\Omega^{1,0}\tens_A F$ with 
$\del_F(\bar\del v)=0$. As $v^*\in V$ for $v\in V$ the formula for the inner product shows that $g_F$ is a hermitian matrix with complex entries. Then $\bar\del g=0$ so all the Christoffel symbols of $\nabla_F$ vanish. For the zero torsion, we apply the formula for torsion to $\extd v$. 

Part (2) is given by (\ref{psidag}). When both parts (1) and (2) apply,  $\nabla_E-\sigma_E\,\dag\,\nabla_F\,*$ is a difference of two connections on $E$, and thus must be a left module map. However it vanishes on all $\del v$ and so must vanish on all of $E$. Then we note that this is the usual definition of $*$-preserving.  \endproof

Now we consider bilinear (not hermitian) `round bracket' metrics defined from a hermitian one by $\<\xi,\overline{\eta}\>=(\xi,\eta^*)$ (e.g. given as in (\ref{eq47}) directly by $\phi_\pm$). First suppose
for $F=\Omega^{0,1}$ that we have a Hermitian metric preserving connection $\nabla_F:F\to\Omega^1\tens_A F$ for $\<\ ,\ \>_F:F\tens \overline{F}\to A$. Then the right connection $\nabla^R{}_E:E\to E\tens_A \Omega^1$ defined by
$\nabla^R{}_E=\dag\, \nabla_F\, *$ together with $\nabla_F$ preserve $(\,,):F\tens E\to A$ in the sense that
\begin{align}
(\id \tens (\,,))  (\nabla_F\tens\id)    + ((\,,)\tens\id)(\id\tens \nabla^R{}_E)=\extd (\,,):F\tens E\to \Omega^1.
\end{align}
Note that this does not require that $(\,,)$ descends to be defined on $\tens_A$. Now suppose that in addition we have the above construction with $E$ and $F$ swapped. There are now left and right connections on $E$ and $F$, and it is reasonable to ask if they are related. In the case where the connections $\nabla_E$ and $\nabla_F$ are bimodule connections, if equation (\ref{tydo})  holds (without any of the assumptions in Proposition~\ref{propBdouble} itself) then we have 
 the following relation between $\nabla_E$ and $\nabla_F$:
\begin{align} \label{starer}
\sigma_E^{-1}\nabla_E=\dag\, \nabla_F\, *,\quad \sigma_F^{-1}\nabla_F=\dag\, \nabla_E\, *.
\end{align}
This says that the right connections $\nabla^R{}_E$ and $\nabla^R{}_F$ for preservation of $(\,,)$ are simply the usual right connection associated to a left bimodule connection. 
Now we put these together by defining a left connection $\nabla:\Omega^1\to \Omega^1\tens\Omega^1$ on
$\Omega^1=E\oplus F$ by $\nabla=\nabla_E\oplus\nabla_F$. This is a bimodule connection if $\nabla_E$ and $\nabla_F$ are, with
$\sigma=\sigma_E\oplus\sigma_F$. Putting this together with (\ref{starer}) gives
\begin{align} \label{starer2}
\sigma\,\dag\, \sigma=\dag,\quad
\sigma^{-1}\nabla=\dag\, \nabla\, *.
\end{align}
These are the usual equations defining a $*$-preserving connection on the $*$-calculus $\Omega^1$. 
Moreover, if we define $(\,,):\Omega^1\tens\Omega^1\to A$ using the decomposition $\Omega^1=E\oplus F$
by $(\,,):F\tens E\to A$ and $(\,,):E\tens F\to A$ and zero on $E\tens E$ and $F\tens F$. Then the bimodule connection $\nabla$ above preserves the metric, i.e.\
\begin{align}  \label{ode}
(\id \tens (\,,))  (\nabla\tens\id)    + ((\,,)\tens\id)(\id\tens \sigma^{-1}\nabla)=\extd (\,,):\Omega^1\tens \Omega^1\to \Omega^1.
\end{align}
Note that this does not require that $(\,,)$ descends to be defined on $\tens_A$, but if it does it gives the usual metric preserving condition.

\subsection{Metric compatible connection on  quantum-braided plane}\label{secexplane}
Recall that for $\C_q^2=B:=B_+(V)$ we used a basis $z$ and $w=q^{-1/2}z^*$ for $V$. We have a comodule star map $\medstar:V\to \overline{V}$
(defined by $v\mapsto \overline{v^*}$) given by
\[
z\mapsto q^{1/2}  \,\overline{w}, \quad w\mapsto     q^{-1/2} \,\overline{z},
\]
where the open font star and overlines are for book-keeping so as to be able to work with bimodule maps. By definition $\Omega^{0,1}=B\tens V$ and $\Omega^{1,0}=V\tens B$ , and then
\[
\bar\partial z=1\tens z,\quad \bar\partial w=1\tens w,\quad \partial z=z\tens 1,\quad \partial w=w\tens 1.
\]
We use $\eta$ the $q$-epsilon tensor
\begin{align} \label{nott}
\eta(z\tens z)=\eta(w\tens w)=0,\quad \eta(z\tens w)=q,\quad \eta(w\tens z)=-1   
\end{align}
which has its image in the centre and  is a morphism for coaction of $\C_q[SL_2]$ (and hence $\C_q[SU_{1,1}]$, but not a morphism for the coaction of the central extension $H$). The braiding, however, is not involutive and the resulting inner product is not expected to (and indeed does not) descend to $\tens_{\C_q^2}$. The map $\eta$ also obeys the condition $*\circ\eta=\eta\circ\dagger$  to have a hermitian inner product from calculate $\phi_-(1\tens u\tens v\tens 1)=\eta(u\tens v)$ in (\ref{eq47}), finding in our case
\begin{align} \label{notuu}
\<\bar\partial z,\overline{\bar\partial z}\>&= \phi_-\big(\bar\partial z\tens\partial (z^*)\big)=q^{1/2}\, \phi_-\big(\bar\partial z\tens\partial w\big)=
q^{3/2},\cr
\<\bar\partial z^*,\overline{\bar\partial z}\>&=\<\bar\partial z,\overline{\bar\partial z^*}\>=0,\quad
\<\bar\partial z^*,\overline{\bar\partial z^*}\> = -q^{1/2}.
\end{align}
Thus on $F=\Omega^{0,1}$ using the basis order $(\bar\partial z,\bar\partial z^*)$ we have
\[
g_F=\left(\begin{array}{cc}q^{3/2} & 0 \\0 & -q^{1/2}\end{array}\right),
\]
which has an indefinite signature. We now carry out the construction from Corollary~\ref{anti-cherncor}. As $\del\bar\del z=\del\bar\del z^*=0$ we have $\del_F\bar\del z=\del_F \bar\del z^*=0$. Then as $\bar\partial g_F=0$ we have 
$\bar\del_F\bar\del z=\bar\del_F \bar\del z^*=0$, so the Chern connection $\nabla_F=\del_F+\bar\del_F$ vanishes on the basis $(\bar\partial z,\bar\partial z^*)$. Evaluating on this basis shows that the torsion of $\nabla_F$ vanishes.

Next, using $\Theta^{1,1}$ as given in Corollary~\ref{corqplane} to compute $\phi_+$ using (\ref{eq47}), we obtain the hermitian metric on $E=\Omega^{1,0}$ as
\[
\< \del z^*,\overline{ \del z^*}\>=q^{1/2},\quad \< \del z,\overline{ \del z}\>= - q^{ 3/2},\quad 
 \< \del z^*,\overline{ \del z}\>= \< \del z,\overline{ \del z^*}\>=0.
\]
Proceeding as before, we get a zero torsion Chern connection $\nabla_E=\del_E+\bar\del_E$ which vanishes on the basis $(\partial z,\partial z^*)$.
Combining these we take a hermitian metric on $\Omega^1=\Omega^{0,1}\oplus \Omega^{1,0}$ by defining $\Omega^{0,1}$ and $\Omega^{1,0}$ to be perpendicular, and there is a resulting zero torsion metric preserving connection $\nabla=\nabla_F\oplus\nabla_E$. Here, $\nabla=0$ on all the basis elements $\{\bar\partial z,\bar\partial z^*,\partial z,\partial z^*\}$. This is a special case of the construction in Proposition~\ref{propBdouble}. Moreover, we can view $\phi_\pm$ as defining a metric $(\, , )$ on $\Omega^1$, and by the discussion after 
Proposition~\ref{propBdouble} $\nabla$ is a $(\, , )$ preserving $*$-preserving bimodule connection, as defined in (\ref{ode}). 
Note that in this case neither $\<\,,\>$ nor $(\, , )$ descend to be defined on $\tens_{\C_q^2}$.

\subsection{Degenerate example on quantum-braided plane that descends}

As a different choice of inner product in Proposition~\ref{propphim}, we set $\eta:V\tens V\to V\tens V\in B=B_+(V)$ to be the quotient map  from the tensor algebra to $B_+(V)$. This is applies in general, but we look at it in the quantum plane case $B=\C_q^2$. Recalling that in the category $\bar\del v=1\tens v\in B\tens V$ and $\del v=v\tens 1\in V\tens B$ this translates to
\begin{align*}
\phi_-(\bar\del z   \tens  \del  z) &= z^2,\quad \phi_-(\bar\del z^*   \tens  \del  z^*) =  z^{*2},\quad \phi_-(\bar\del z   \tens  \del  z^*) = z\, z^*,\quad \phi_-(\bar\del z^*   \tens  \del  z) =  z^* z.
\end{align*}
The reader may check that this gives a bimodule map $\phi_-:\Omega^{0,1}\tens_B \Omega^{1,0}\to B$. 
From a formula in the proof of Proposition~\ref{propphim} we have $\phi_+(v\tens 1\tens 1\tens v')=\eta\circ\Psi^{-1}(v\tens v')$, which gives, using the formula for the braiding in the proof of Lemma~\ref{lemqplane*} and recalling that $z^*=q^{{1\over 2}}w$,
\begin{align*}
\phi_+(\del z   \tens \bar \del  z) &=  q^{-2} z^2,\quad \phi_+(\del z^*   \tens  \bar \del  z^*) =  q^{-2} z^{*2} ,\\
  \phi_+(\del z   \tens \bar \del  z^*) &= 
q^{-2} z\,z^* ,\quad \phi_+(\del z^*   \tens \bar \del  z) =  q^{-2} z^* z .
\end{align*}
On $\Omega^{1,0}$, if we set $e^1=\del z$ and $e^2=\del z^*$ then we have the hermitian matrix
\[
(g^{ij}) = q^{-2} \left(\begin{array}{cc}z\,z^* & z^2 \\z^{*2} & z^* z\end{array}\right).
\]
As $e^1$ and $e^2$ are free generators of $\Omega^{1,0}$ as a left module, the idempotent $P$ is the identity matrix, so 
$(\tilde g_{ij})$ would have to be the matrix inverse of $(g^{ij})$. However $(g^{ij})$
 is not invertible, even if we adjoin inverses $z^{-1}$ and $z^{*-1}$, so the inner product is degenerate. 
 This is the unique metric, up to scaling, which descends to $\Omega^{0,1}\tens_B \Omega^{1,0}$ which takes values in polynomials in $z$ and $z^*$. For a  connection, we cannot use our formulae for the Chern connections without the inverse metric but one can, in principle, directly solve for connections on $\Omega^{1,0}$ and $\Omega^{0,1}$ with similar properties without expecting these to be unique.

\subsection{The doubled construction in the finite group case}

We begin with the conjugate version of Proposition~\ref{propchernG}.

\begin{lemma}\label{lemchernG} In the context of Proposition~\ref{propGdolb} with both $\CC^{1,0},\CC^{0,1}$ conjugacy classes and (\ref{hur}) imposed, suppose that we are given a nondegenerate hermitian metric of the form $\<e^{a}, e^{b}\>=g^{ab}\in C(G)$ on $F=\Omega^{0,1}$, where $a,b\in \CC^{0,1}$. Here, $F$ is spanned over $A=C(G)$ by $\{e^{a}\}$ for $a\in \CC^{0,1}$. Then we have a hermitian metric preserving connection on $F$ by
 \[
\nabla_F e^a=\sum_{b,c\in\CC^{0,1}} \bar\partial g^{ab}\ \tilde g_{bc} \tens e^c,
 \]
with torsion
  \[
T_F e^a=\sum_{b,c\in\CC^{0,1}} \bar\partial g^{ab}\ \tilde g_{bc} \wedge e^c - \bar\partial e^a.
 \]
(Note that in this case $\tilde g_{bc}$ is the matrix inverse of $g^{ab}$.) Moreover, the hermitian metric descends to $\tens_A$ if and only if 
the matrix $g^{ab}$ is diagonal, and in this case we have a bimodule connection with
 \begin{align*}
\sigma_F(e^a\tens e^{b^{-1}})=e^{ab^{-1}a^{-1}}\tens e^a,\quad \sigma_F(e^a\tens e^b)=\frac{g^{aa} } {R_{aba^{-1}}(g^{aa})}\, 
e^{aba^{-1}}\tens e^a.
\end{align*}
\end{lemma}
\proof Begin with, for $a\in  \CC^{0,1}$, and use that $\CC^{1,0}$ is a union of conjugacy classes:
\[
\del e^a=\sum_{b\in  \CC^{0,1}} \big(e^a\wedge e^{b^{-1}} + e^{b^{-1}} \wedge e^a\big) = 
\sum_{b\in  \CC^{0,1}} \big(e^a\wedge e^{b^{-1}} + e^{ab^{-1}a^{-1}} \wedge e^a\big)=0
\]
due to (\ref{hur}). (The asymmetry introduced by (\ref{hur}) is the reason that this result is simpler than Proposition~\ref{propchernG}.)
Then, following the method of Corollary~\ref{cherncor}, we have 
 $\partial_F e^a=0$. Thus $\Gamma_-=0$ and $\Gamma_+=-(\bar\partial g)\,\tilde g$, so $\bar\del_E e^a = \nabla_E e^a$ gives the formula in the statement for a metric preserving connection. Here, $\tilde g$ is the matrix inverse as the algebra is commutative and the generators $e^a$ are free. Also as the algebra is commutative, if the metric descends to $\tens_A$ then we can only have nonzero values of 
 $\<e^a,\overline{e^b}\>$ if $ab^{-1}=e$, so in this case $g^{ab}$ is diagonal. Then, for $f\in A$ and $a\in \CC^{0,1}$, we consider $\sigma_F(e^a\tens\extd f)=\nabla_F(e^a\,f)-\nabla_F(e^a)\, f$ to get
  \begin{align*}
\sum_{c\in\CC} &(R_{ac}(f)-R_a(f))\,\sigma_F(e^a\tens e^c) = \extd(R_a(f))\tens e^a
+ R_a(f)\,\nabla_F(e^a) - \nabla_F(e^a) \,f \cr
&= \sum_{c\in\CC}  (R_{ca}(f)-R_a(f))\, e^c\tens e^a + R_a(f)\,\bar\partial g^{aa}\ \tilde g_{aa} \tens e^a
- \bar\partial g^{aa}\ \tilde g_{aa} \tens e^a \, f\cr
&= \Big(\sum_{c\in\CC}  (R_{ac}(f)-R_a(f))\, e^{aca^{-1}} + R_a(f)\,\bar\partial g^{aa}\ \tilde g_{aa} 
- \bar\partial g^{aa}\ \tilde g_{aa}\, R_a(f) \Big) \tens e^a .
\end{align*}
For $c\in\CC^{1,0}$, we have the first part of $\sigma_F$ as stated. If $c=b\in \CC^{0,1}$ then we have
  \begin{align*}
&\sum_{b\in\CC^{0,1}}   \Big((R_{ab}(f)-R_a(f))\, e^{aba^{-1}} + \frac{ R_b(g^{aa})-g^{aa}}{R_b(g^{aa})} \, (R_a(f)\,e^b- e^b\, R_a(f))
\Big) \tens e^a \cr
&=\sum_{b\in\CC^{0,1}}   \Big((R_{ab}(f)-R_a(f))\, e^{aba^{-1}} +\Big( 1 - 
 \frac{ g^{aa}}{R_b(g^{aa})} \Big)
 \, (R_a(f)-R_{ba}(f))\,e^b      \Big)    \tens e^a \cr
&=\sum_{b\in\CC^{0,1}}   \Big((R_{ab}(f)-R_a(f))\, e^{aba^{-1}} +
\Big( 1 - 
 \frac{ g^{aa}}{R_{aba^{-1}}  (g^{aa})} \Big)
  \, (R_a(f)-R_{ab}(f))\,e^{aba^{-1}}
\Big) \tens e^a,
\end{align*}
giving the second part of $\sigma_F$. 
 \endproof

Similarly to Corollary~\ref{doubleB}, we can put the $E=\Omega^{1,0}$ and $F=\Omega^{0,1}$ results together in the group case.

\begin{corollary} In the context of Proposition~\ref{propchernG} for the algebra $C(G)$, we take subsets $\CC^{1,0}$ defining the calculus
$E=\Omega^{1,0}$ and $\CC^{0,1}$ defining $F=\Omega^{0,1}$ with hermitian metrics $g_E,g_F$ which we view on  $\Omega^1=E\oplus F$ with $E$ and $F$ perpendicular. The left connection $\nabla=\nabla_{E\oplus F}$ on $\Omega^1$ restricting to $\nabla_E$ and $\nabla_F$ on the two parts is hermitian metric compatible for $\<\, , \>$ and has torsion restricting to $T_E$ and $T_F$. If both metrics $g_E$ and $g_F$ are diagonal then $\nabla$ is a bimodule connection, with corresponding $\sigma$ given by the sum of $\sigma_E$ and $\sigma_F$. 
\end{corollary}
\proof  The important thing distinguishing $\nabla=\nabla_{E\oplus F}$  from a general connection on $\Omega^1$ is that for $e\in E$ we have $\nabla e\in \Omega^1\tens_A E$ and for $f\in F$ we have $\nabla f\in \Omega^1\tens_A F$. This, together with $\<\overline{e},f\>=\<\overline{f},e\>=0$, is enough to show that $\nabla$ is metric preserving by considering the 
cases for $\Omega^1=E\oplus F$. Similarly, we consider the cases for the inputs to the torsion and $\sigma$. 
\endproof

 Note that in the $C(G)$ case with diagonal metrics, the diagonal entries define arrow weights and it is natural to ask for the metric to be to edge-symmetric\cite{BegMa} so that these do not depend on the arrow direction. Here $g^{ab}=-h^{ab^{-1}}=\delta_{a,b}c_a$, say, corresponds to  $h^{ab}=(e^a,e^b)$, and edge symmetry is then $c_{a^{-1}}=R_a(c_a)$. This is for $a,b\in \CC$ on any finite group. In our case it means for $a,b\in C^{0,1}$ that we need 
 \[ g^{a^{-1}b^{-1}}=R_a (g^{ab})\]
 as the relation between the two metrics.

\begin{example}\label{exZ} We apply the theory to $G=\Z$ and $\CC=\{1,-1\}$, taking $E=\Omega^{1,0}$ to have generator $e^+$ and 
$F=\Omega^{0,1}$ to have generator $e^-$. The relations for $\Omega_{LL}$ are $e^{\pm} \wedge e^{\pm} =0$ and
$e^{+} \wedge e^{-} + e^{-} \wedge e^{+} =0$, so we get $\extd e^{\pm}=0$. (In this case we automatically get condition (\ref{hur}).)
From this
$\bar\del_E(e^+)=0$ and $\del_F(e^-)=0$. The hermitian metric on $E$ is given by a real valued function $g^{++}$ so that
$\<e^+,\overline{e^+}\>_E=g^{++}$. From 
Corollary~\ref{cherncor} the matrix of Christoffel symbols of $ \partial_E$ is
\[
\Gamma_+= - \partial(g)\,\tilde g = -(\partial g^{++}(g^{++})^{-1})
\]
so
\[
\nabla_E e^+=\del_E e^+=\big(1-\frac{ g^{++}  }{    R_{+1}(g^{++})   }\big)\, e^+ \tens e^+.
\]
The metric matrix is trivially diagonal (as it is $1 \times 1$), so from 
Proposition~\ref{propchernG} we have a bimodule connection and an expression for $\sigma_E$. 
Similarly, the hermitian metric on $F$ is given by a real valued function $g^{--}$ so that
$\<e^-,\overline{e^-}\>_F=g^{--}$, and the matrix of Christoffel symbols of $\bar \partial_F$ is
\[
\Gamma_-= - \bar\partial(g)\,\tilde g = -(\bar\partial g^{--}(g^{--})^{-1})
\]
so
\[
\nabla_F e^-=\bar\del_F e^-=\big(1-\frac{ g^{--}  }{    R_{-1}(g^{--})   }\big)\, e^- \tens e^-.
\]
Hence, we obtain $\nabla=\nabla_E\oplus\nabla_F$ as
\[
\nabla(e^{\pm})=(1-\rho_{\pm})\,  e^{\pm}\tens e^{\pm},\quad \rho_\pm=\frac{ g^{\pm\pm}}{    R_{\pm1}(g^{\pm\pm})}
\]
and can calculate
\[
\sigma(e^{\pm}\tens e^{\mp}) = e^{\mp}\tens e^{\pm},\quad 
\sigma(e^{\pm}\tens e^{\pm}) = \rho_\pm  e^{\pm}\tens e^{\pm}.
\]
This is hermitian metric compatible by construction and one can check that it is $*$-preserving, hence also metric compatible for $(\ ,\ )$. Indeed, we obtain exactly the QLC found in \cite{Ma:haw} where $\cg=-g_+e^+\tens e^-- g_-e^-\tens e^+$ with $R_{\pm 1}(g_\mp)=g_\pm$ as the metric is edge symmetric, and  $g^{\pm\pm}=1/R_{\pm 1}(g_\mp)$. 
\end{example}

\section{Inner products with values in bimodules}\label{secbundle}

In this section, we introduce the idea of inner products that have values in an $A$-bimodule $L$ rather than in the coordinate algebra $A$. We first make some general remarks in the hermitian case then focus on a specific construction over any field. 

\subsection{Hermitian inner products with values in line bundles} \label{sechermbundle}
Let $A$ be a $*$-algebra with a $*$-differential calculus, and $L$ be a line bundle with conjugate $\overline{L}$ the same as its dual. This means we have isomorphisms $\mu:L\tens_A \overline{L}\to A$ and $\nu: \overline{L}\tens_A L\to A$ so that $\id.\mu=\nu.\id:
\overline{L} \tens_AL\tens_A \overline{L} \to \overline{L}$ and $\mu.\id=\id.\nu:
L\tens_A \overline{L} \tens_A L\to L$. We also assume that we have an involutive star operation on $L$, and this gives a bimodule isomorphism $\medstar:L\to \overline{L}$ by $\medstar(\ell)= \overline{\ell^*}$. There is an antilinear star operation on $\overline{L}$
given by $(\overline{\ell})^*=\overline{\ell^*}$. Finally we assume that $\mu$ and $\nu$ are connected by the star operation, i.e.\ 
$\mu(m\tens\overline{\ell})^*=\nu(\overline{\ell^*} \tens m^*)$. This means that the tensor algebra of $L$ and $\overline{L}$ can be made into a graded star algebra \cite{BBLine}, where we use the associativity rule above for $\mu$ and $\nu$ to replace these operations by a product.

Let $E$ be an f.g.p.\  left $A$-module with dual basis $e_i\tens e^i\in {}_A\Hom(E,A)\tens E$. We write the evaluation map
as $\ev:E\tens {}_A\Hom(E,A)\to A$.
We then have $e=\ev(e\tens e_i).e^i$ for all $e\in E$ and $f=e_i.\ev(e^i\tens f)$ for all $f\in  {}_A\Hom(E,A)$. 
Then $\ev(e^i\tens e_j)=P^i{}_j$ is a matrix valued idempotent $P\in M_n(A)$ (supposing that there are $n$ elements in the dual basis). Using a right module map
\[
G:\overline{E}\to {}_A\Hom(E,A)\tens_A L
\]
we define an inner product as an $A$-bimodule map (where $e,e'\in E$)
\begin{align} \label{urg}
\<\, ,\>: E\tens \overline{E} \to L,\quad \<e',\overline{e}\>= (\ev.\id)(e'\tens G(\overline{e})).
\end{align}
As $E$ is f.g.p.\  the map $G$ must have the form $G(\overline{e^i})=e_j\tens g^{ji}$ for some $g^{ji}\in L$, so we get a matrix $g\in M_n(L)$. 
As, in $ {}_A\Hom(E,A)\tens_A L$,
\[
e_j\tens g^{ji}= e_k\, P^k{}_j\tens g^{ji}= e_k\tens P^k{}_j\, g^{ji}
\]
we assume, without loss of generality, that $P\,g=g\in M_n(L)$. Next,
\[
e_j\tens g^{ji}=G(\overline{e^i})=G(\overline{P^i{}_j\,e^j})=G(\overline{e^j})\, (P^i{}_j)^*= e_k\tens g^{kj}\, (P^i{}_j)^*
\]
and evaluating on an element $e^q$ on the left we get $P\,g=P\,g\,P^*$, where $P^*\in M_n(A)$ is the star (on individual entries) transpose matrix, i.e.\ the usual matrix star in the case of ordinary matrices. Thus we deduce that we must have $g=g\,P^*$, and this is enough to ensure that the formula $G(\overline{e^i})=e_j\tens g^{ji}$ extends to a right module map by
\[
G(\overline{e})= G(\overline{\ev(e\tens e_i).e^i}) = G(\overline{e^i})\, \ev(e\tens e_i)^*.
\]

The inner product (\ref{urg}) will be called hermitian when $\<e,\overline{e'}\>=\<e',\overline{e}\>^*$. As 
\begin{align}
\<e^i,\overline{e^j}\> = \ev(e^i\tens e_k) g^{kj} = P^i{}_k\, g^{kj} = g^{ij},
\end{align}
so for a hermitian inner product we require $(g^{ij})^*=g^{ji}$, or in terms of matrices $g^*=g\in M_n(L)$. The inner product (\ref{urg}) will be called nondegenerate when $G$ is invertible, and we write (using some manipulation of line bundles) the inverse as $H:
{}_A\Hom(E,A) \to \overline{E}\tens_A \overline{L}$, where
\[
(\id.\nu)(H\tens\id)G(\overline{e}) = \overline{e},\quad (\id.\mu) (G\tens\id)H(f) = f.
\]
We write
\[
H(e_i)=\overline{e^j}\tens \tilde g_{ji},
\]
where $\tilde g\in M_n(\overline{L})$. As before we can assume that $P^*\,\tilde g=\tilde g$ and then deduce that $\tilde g\,P=\tilde g$. 
Then the inverse definition gives
\[
e_i= (\id.\mu) (G\tens\id)H(e_i) = G(\overline{e^j}  )   . \tilde g_{ji} = e_k.(g^{kj}. \tilde g_{ji}).
\]
 Evaluating this expression against $e^q$ on the left gives $P=P\,g\,\tilde g$, so we deduce that $g\,\tilde g=P$.  (This corresponds in classical geometry to $\tilde g$ being the inverse matrix to $g$.) Then we have the inverse the other way round,
\[
\overline{e^i}=(\id.\nu)(H\tens\id)G(\overline{e^i}) = H(e_j). g^{ji} = \overline{e^k} . \tilde g_{kj}.g^{ji} = 
 \overline{   (   \tilde g_{kj}.g^{ji} )^*\,   e^k} 
\]
so we have $e^i= (   g^{ji} )^*\,   (   \tilde g_{kj})^* e^k$ and evaluating this on $e^q$ on the right gives
\[
P=g^*\,\tilde g^*\,P= g^*\, (P^*\,\tilde g)^*= g^*\, \tilde g^*= (\tilde g\, g)^*
\]
so we deduce that $\tilde g\, g=P^*$. Finally
\[
\tilde g = P^*\, \tilde g = \tilde g^*\, g^* \tilde g=\tilde g^*\, g\, \tilde g = \tilde g^*\, P=(P^*\, \tilde g)^* =\tilde g^*.
\]

For an example, we will shortly construct an $L$-valued metric on $\Omega^1$ for the quantum plane with its 4D calculus which, via $*$, can be viewed as an $L$-valued hermitian metric on $\Omega^1$. 

\subsection{Connections preserving line bundle valued metrics}
Suppose that we have a left bimodule connection $(\nabla_L,\sigma_L)$ on the line bundle $L$. Then we can give a condition for a left connection $\nabla_E:E\to \Omega^1\tens_A E$ to preserve an $L$-valued hermitian metric as the vanishing of
\[
M=(\id\tens\<\,,\>))(\nabla_E\tens\id) + \sigma_L(\<\,,\>\tens\id)(\id\tens\nabla_{\overline{E}})-\nabla_L \<\,,\>:E\tens\overline{E}\to
\Omega^1\tens_A L.
\]
The operator  $M$ is a bimodule map, we check the more difficult direction.
\begin{align*}
M(e'\tens \overline{e}\, a)- M(e'\tens \overline{e})\, a = \sigma_L(\<e', \overline{e}\>\tens\extd a) - \nabla_L(
\<e', \overline{e}\>\, a) + \nabla_L(\<e', \overline{e}\>)\, a,
\end{align*}
and this vanishes by definition of $\sigma_L$. As $M$ is a bimodule map, it is enough to check that it vanishes on the dual basis. As before we set $\nabla_E e^i= - \Gamma^i{}_j\tens e^j$, and then
\begin{align*}
M(e^i\tens\overline{e^j}) &= - \Gamma^i{}_k\tens g^{kj} - \sigma_L\big( g^{ik}\tens (\Gamma^j{}_k)^*\big) -\nabla_L(g^{ij}).
\end{align*}
The vanishing of this is equivalent to hermitian metric preservation.

\subsection{Inner product $(\, , )$ on $\Omega^{1,0}$}

We can also use similar methods as in Proposition~\ref{propphim} to construct a natural inner product on $\Omega^{1,0}$, for $A=B:=B_+(V)$, e.g. for the standard quantum plane 2D calculus, and working over a field (we do not necessarily have to use a $*$-calculus or to double up). We assume a map $\phi_-:\Omega^{0,1}\tens\Omega^{1,0}\to B$ built from $\eta:V\tens V\to B$ as in (\ref{phim}) but we do not require that this itself descends to $\tens_B$. 

\begin{figure}
\[ \includegraphics[scale=0.9]{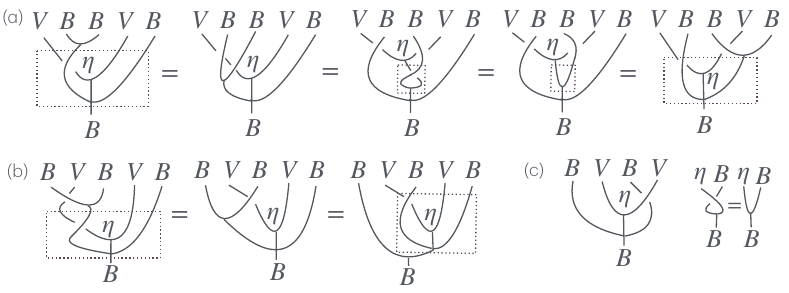}\]
\caption{Proof of  Proposition~\ref{leftmetric}. Here $B$ stands for $B_+(V)$. \label{fig4}.}
\end{figure}

\begin{proposition}\label{leftmetric} For $B:=B_+(V)$, if $\eta:V\tens V\to B$ is a morphism  then $\Omega^{1,0}=V\und{\tens} B$ acquires a bimodule map 
\[ \psi_+=\phi_-\circ(\Psi_{B,\Lambda^1}^{-1}\tens\id):\Omega^{1,0}\tens\Omega^{1,0}\to B\]
 and descends to $\tens_B$ if the image of $\eta$ is braided-commutative as in Figure~\ref{fig4}(a).  Similarly for $\psi_-$ on $\Omega^{0,1}=B\tens V$ in Figure~\ref{fig4}(c).
\end{proposition}
\proof The proof that $\psi_+$ is a right module map is immediate without any conditions on $\eta$. Likewise that it is a left module needs no conditions on $\eta$ and is shown in Figure~\ref{fig4}(b). We use that the product of $B$ is a morphism and associativity. So we always have a map $\Omega^{1,0}\tens\Omega^{1,0}\to B$. Part (a) shows that this descends to $\tens_B$ using that the product is a morphism and then that $\eta$ is a morphism (as well as associativity steps omitted). The third equality is braided-commutativity in the sense of the two small boxed expressions being equal. We then use associativity again. Finally, Figure~\ref{fig4}(c) gives the parallel right-handed version for 
\[ \psi_-=\phi_-\circ(\id\tens\Psi_{\Lambda^1,B}^{-1}):\Omega^{0,1}\tens_B\Omega^{0,1}\to B\]
 and the appropriate braided-commutativity condition on the image of $\eta$ for this to descend to $\tens_B$ (the proofs are reflected followed by a reversal of all braid crossings). That $\psi_\pm$ are bimodule maps also follows from the bimodule map property of $\phi_-$ and left/right module map properties of the relevant inverse braidings. \endproof

Although we have used a double-complex notation, the result applies to each half separately. However, if we do assume that $V$ is unitary and work in the $\Omega^{p,q}$ factorisable double complex with $B_+(V)$ a $*$-algebra then $(\, , )$ defined as $\psi_\pm$ on the relevant components and zero elsewhere is a `real' generalised metric in the sense of \cite{BegMa} if and only if 
\[ *\circ \psi_+=\psi_-\circ\dagger,\]
which using (\ref{psidag}) holds if and only if  $*\circ\phi_-=\phi_-\circ\dagger$ as in Proposition~\ref{propphim} for the condition on $\eta$ given there. Moreover, if  the braided-commutativity in Proposition~\ref{leftmetric} holds then $(\, , )$ descends to $\tens_B$. It is not assumed to be quantum symmetric or invertible. 

Next, as with Proposition~\ref{propphim}, the diagrammatic proof of Proposition~\ref{leftmetric} works in exactly the same way  if we let  $\eta: V\tens V\to L$ where $\eta$ is a morphism in the category and $L$ is a bimodule in the category, but this time we do not require that $\Psi$ is involutive. We still need that the image of $\eta$ is braided-commutative in the same way for the 3rd equality in Figure~\ref{fig4} (a) and then obtain $\psi_+:\Omega^{1,0}\tens_B\Omega^{1,0}\to L$. Similarly, if the image of $\eta$ obeys the different braided-commutativity in Figure~\ref{fig4} (c) then we obtain $\psi_-:\Omega^{0,1}\tens_B\Omega^{0,1}\to L$.

Now suppose that the category we are working in is the comodules of a central extension $H=G\tens\C_\alpha[t,t^{-1}]$ with coquasitriangular structure extended by $\CR(t^m,t^n)=q^{\alpha mn}$. At the moment we are not considering the bar category structure, this works over any field but we have written $\C$ as the case of interest. Now let $\eta:V\tens V\to \C$ be invariant under $G$, where $G$ right coacts on $V$ and $H$ right coacts with an extra $\tens t$. Let $\C D$ be a 1-dimensional comodule with $\Delta_R D=D\tens 1\tens t^2$ and tensor $\eta$ by $D$ to a map which we again denote $\eta: V\tens V\to \C D$. Now $\eta$ is a morphism in the category. 

\begin{corollary}\label{corLpm}Let $B:=B_+(V)$ and $\eta:V\tens V\to \C D$ be a morphism in the category as above. Then
 \[ L_\pm=B\tens \C D,\quad  b.(c\tens D)=bc\tens D,\quad  (c\tens D).b=c\Psi^{\mp 1}(D\tens b)\]
 are bimodules for which the image of $\eta$ is braided-commutative in the required way and  give
 \[ \psi_+: \Omega^{1,0}\tens_B\Omega^{1,0}\to L_+,\quad  \psi_-: \Omega^{0,1}\tens_B\Omega^{0,1}\to L_-.\]
 \end{corollary}
 \proof The proofs that we obtain bimodules with the appropriate braided-commutativity of the image of $\eta$ are easily done by diagrams and left to the reader. Note that that the stated right actions for the bimodule structures are explicitly $(b\tens D).c= q^{\mp 2\alpha |c|} bc\tens D$ where $|c|$ is the degree of $c\in B_+(V)$, from which one can also verify everything.  We then apply the bimodule versions of Proposition~\ref{leftmetric} for $\psi_\pm$. 
 \endproof
 
It is also convenient to let $\tilde B=B\rtimes\C Z=B[\delta,\delta^{-1}]$ where we adjoin an invertible element $\delta$ to $B=B_+(V)$ with commutation relations $v\delta=q^{2\alpha}\delta v$. Then we can view $L_\pm$ as submodules of grade $\pm 1$ in $\tilde B$, where the grading is the by powers of $\delta$. Here
\[ i_+(b \tens D)=b\delta,\quad i_-(b\tens D)=b\delta^{-1},\quad i_\pm: L_\pm\hookrightarrow \tilde B\]
which one can check are bimodule maps when we view $B\subset\tilde B$ and multiply there. In this way, we can view both $\psi_\pm$ as inner products valued in $\tilde B$. In the bar category setting one can relate the two inner products via $*$. We illustrate this, and the construction of a Levi-Civita connection, in the case of  the quantum-braided plane.

\subsection{Example of $(\, , )$ on $\Omega^{1,0}$ on the quantum plane}\label{secexplaneround} We use the $\C_q[SU_{1,1}]$ invariant $\eta:V\tens V\to \C$ as in (\ref{nott}), which is not fully covariant under the centrally extended quantum group, i.e. we are in the setting of the Corollary~\ref{corLpm}. This  gives us a  bimodule map $\psi_+:\Omega^{1,0} \tens_B \Omega^{1,0}  \to \tilde B$ with $\alpha=3/2$ so that $v\delta=q^3\delta v$ for all $v\in V$. The result is given on the generators as 
\begin{align*}
& \psi_+(\del z^*\tens\del z)=- q^{1/2}\,\delta,\quad \psi_+(\del z\tens\del z^*)=q^{3/2}\, \delta,\quad  \psi_+(\del z^*\tens\del z^*)=\psi_+(\del z\tens\del z)=0. 
\end{align*}
We also have the bimodule map $\psi_-:\Omega^{0,1} \tens_B \Omega^{0,1}  \to B\,\delta^{-1}$ given by
\begin{align*}
& \psi_-(\bar\del z^*\tens\bar\del z)=- q^{1/2}\,\delta^{-1},\quad \psi_-(\bar\del z\tens\bar\del z^*)=q^{3/2}\, \delta^{-1},\quad \psi_-(\bar\del z^*\tens\bar\del z^*)=\psi_-(\bar\del z\tens\bar\del z)=0.
\end{align*}
As an useful check, one can see that these maps indeed descend to $\tens_B$, for example
\begin{align*}
& \psi_-(\bar\del z\, z^*\tens\bar\del z^*) = q\,  \psi_-(z^*\,\bar\del z\tens\bar\del z^*) = q\, z^*\,  \psi_-(\bar\del z\tens\bar\del z^*) 
 =q^{5/2}\,z^*\, \delta^{-1},\cr
 & \psi_-(\bar\del z\tens z^*\,\bar\del z^*) = q^{-2}\,  \psi_(\bar\del z\tens\bar\del z^*\, z^*)= q^{-2}\,  \psi_(\bar\del z\tens\bar\del z^*)\, z^*=q^{-1/2}\, \delta^{-1}\,z^*,
\end{align*}
which are equal as $z^* \delta^{-1} = q^{-3}\, \delta^{-1}z^*$. By contrast, this does not work for $\phi_-$ in Section~\ref{secexplane} even if the values stated there are viewed times $\delta$.

Finally we note we can make $\tilde B$ into a $*$-algebra with  $\delta^*=\delta^{-1}$ and in this case we have
\[ \psi_-=*\circ \psi_+\circ\dag\]
In this case we define $(\, , )=\psi_+\oplus \psi_-:\Omega^1\tens_B\Omega^1\to L=B\tens (\C \delta\oplus \C\delta^{-1})\subseteq \tilde B$, where $L\cong L_+\oplus L_-$ as a direct sum bimodule over $B$. Moreover, this is now `real' in the sense $*\circ(\, , )=(\, , )\circ\dagger$ of QRG but with values in the $*$-algebra $\tilde B$ or in $L$ viewed as a bundle with the structure of a star object.

To give a metric preserving connection we begin with a left bimodule connection on $\tilde B$ by defining $\nabla_{\tilde B}(\delta)=\nabla_{\tilde B}(\delta^{-1})=0$. Then 
\[
\sigma_{\tilde B}(\delta^{\pm 1}\tens \extd b) = \nabla_{\tilde B}(\delta^{\pm 1}\, b)- \nabla_{\tilde B}(\delta^{\pm 1})\, b = q^{\mp 3|b|}\, \nabla_{\tilde B}(b\, \delta^{\pm 1})
= q^{\mp 3|b|}\, \extd b \tens \delta^{\pm 1}.
\]
We will implicitly use $\sigma_{\tilde B}$ to swap the orders of $\delta^{\pm}$ and 1-forms by writing $\delta^{\pm}\,\extd b = q^{\mp 3|b|}\, \extd b \ \delta^{\pm 1}$. Now we can
use a left bimodule connection $\nabla$  on $\Omega^1$ to preserve the metric in the sense
\[
(\id .(\,,))(\nabla^L\tens\id) + ( (\,,).\id)(\id\tens\sigma^{-1}\nabla) =\nabla_{\tilde B}\,(\,,):\Omega^1\tens_A\Omega^1 \to \Omega^1.\tilde B.
\]
This equation is analogous to (\ref{ode}), and it is
 satisfied if $\nabla$ vanishes on the basis elements $\del z$, $\del z^*$, $\bar\del z$, $\bar\del z^*$. 
Thus $\nabla$ is exactly the same $*$-preserving bimodule connection as described at the end of
Section~\ref{secdoublech}.

\end{document}